\magnification=1333
\input amstex
\documentstyle{amsppt}
\NoBlackBoxes
\loadbold
\pagewidth{16.5truecm}
\pageheight{23.0truecm}
\NoRunningHeads
\def\curraddr#1\endcurraddr{\address {\it Current address\/}:
#1\endaddress}
\voffset=0.42truecm
\vcorrection{-0.75cm}
\topmatter
\title  \vskip 0.8truecm
            GENERALIZED HERMITE POLYNOMIALS AND
        THE BOSE-LIKE OSCILLATOR CALCULUS        \endtitle
\author            Marvin Rosenblum          \endauthor 
\address           \rm Department of Mathematics  \newline
                   University of Virginia  \newline
                   Math/Astro Building     \newline
                   Charlottesville, Virginia 22903-3199  \newline
                   U. S. A. \newline
                    mr1t\@virginia.edu    \newline \newline  
                    MSC 1991 Primary 33C45, 81Q05 \newline
                   Secondary 44A15 
      \endaddress

\email              mr1t \@ virginia.edu                  \endemail
\dedicatory Dedicated to Moshe Liv\v{s}ic
\enddedicatory
\abstract This paper studies a suitably normalized set of
generalized
Hermite polynomials and sets down a relevant Mehler formula,
Rodrigues
formula, and generalized translation operator. 
 Weighted generalized Hermite polynomials
are the eigenfunctions of a generalized Fourier transform which
satisfies
an  F. and M. Riesz theorem on the absolute continuity  of analytic
measures.
 The Bose-like oscillator calculus, which generalizes the calculus
associated with the quantum mechanical simple harmonic oscillator,
is studied in terms of these polynomials.
\endabstract
will
\endtopmatter
\document
\head 1. Introduction                
\endhead                             
\define\bbH{\Bbb H}
\define\half{\tfrac{1}{2}}
\define\bk#1{\boldkey #1}

\define\pd#1#2{\dfrac{\partial#1}{\partial#2}}
\define\jump{\qquad \qquad \quad \smc }
\define\tild{ \widetilde{\omega}(t) }
      The generalized Hermite polynomials ${\bbH}^{\mu}_n (x),n \in
{\Bbb N} = \{ 0, 1, 2, \dots \},$ 
were defined by Szego \cite{29 \rm, p380, Problem 25} as a set of
real
polynomials orthogonal with respect to the weight $\vert
x\vert^{2\mu}
e^{-x\sp2}$, $\mu > -\half$,  with the degree of $\bbH^{\mu}_n$
equal n. Thus   $\int^ \infty _{-\infty}{\bbH}^ \mu _ m(x)
{\bbH}^\mu _n(x)  e^{-x\sp 2}  \vert x\vert ^{2 \mu} \,dx=0,\ m\neq
n$.
These polynomials can be exhibited in terms of certain
confluent hypergeometric polynomials, or  in terms
of certain generalized Laguerre polynomials.
\endgraf       We refer to Erd\'elyi  \cite{15, Vol 1}
for the definition and properties
of the confluent hypergeometric function $\Phi$ , and  generally
of other special functions. The   m-th confluent hypergeometric
polynomial with parameter $ \gamma + 1   > 0 $ is given by
$$ \align  \Phi (-m, \gamma+1, x) :=& \sum _{k=0}^{m}{ (-1)^k 
\binom mk
 \frac {\Gamma(\gamma + 1)}{\Gamma(k + \gamma  + 1)}
           x^k   } ,\\ \text {which}\  \
 =& \frac {m! \Gamma (\gamma + 1)}{\Gamma (m + \gamma+ 1)}
 L^\gamma_m(x).       \endalign                        $$
Necessarily ${\bbH}_{2m}^\mu (x) = c_{2m} \Phi (-m, \mu + \half,
x^2)$
and ${\bbH}_{2m + 1}^\mu (x) = c_{2m+1} x \Phi (-m, \mu + \frac
{3}{2}, x^2)$
, where $ m \in \Bbb N = \{0,1,2 \dots\},$ and the $ c_.$ are real
constants.
(See Chihara \cite {5,p43,p157}.)
\endgraf In his Ph.D. thesis Chihara \cite{4} normalized these
polynomials
so that the coefficient of $x^n$ in $ {\bbH}_{n}^\mu $  
is $2^n$. Others studying these polynomials,
in general with  varying normalizations, are Dickinson and Warski
\cite {10} and Dutta, Chatterjea, More \cite{13}.  
 We shall set down a different normalization, one that is
appropriate
for our applications.  We shall denote the Chihara polynomials by
$ \{{\bbH}_n^\mu\}_0^\infty $ and our class by
 ${\{H_n^\mu\}_0^\infty} $. \endgraf
We study the generalized Hermite polynomials in section 2. In
section 3
we define a relevant Fourier transform and heat equation, and in
section 4 a relevant translation operator. Finally, in section 5,
we
study some basic aspects of the Bose-like oscillator calculus which
intrinsically  connect with the earlier sections. It will develop
that the Bose-like oscillator calculus is a remarkable, fully
structured generalization of the calculus
associated with the quantum mechanical harmonic oscillator, that
is,
the Boson calculus. Our main result, Theorem 5.12, is a
generalization
of the von Neumann uniqueness theorem to the Bose-like oscillator
calculus.
\head 2. Generalized Hermite Polynomials
\endhead
    Let $ {\Cal C} $ be the set of complex numbers and let
${\Cal C_o} = {\Cal C}\smallsetminus \{-\tfrac 1 2,-\tfrac 3 2,
    -\tfrac 52,\dots  \}. $
\proclaim{\jump    2.1 Definition}  Suppose  $\mu \in {\Cal C_o}$.
 The \bf  generalized Hermite polynomials \sl
${\{H_n^\mu\}_0^\infty} $   are defined for n even by
$$\align H_{2m}^\mu(x)
&:=  (-1)^m \frac {(2m)!}{m!}  \Phi(-m,\mu + \half,x^2) \tag2.1.1 
\\
&= (-1)^m \frac {(2m)!}{m!} \sum _{k=0}^m (-1)^k
          \binom mk \frac {\Gamma(\mu + \half)}{\Gamma(k + \mu +
\half)}
           x^{2k}  .       \endalign             $$
\text {They are defined  for n odd by}
$$ \align   H_{2m+1}^\mu(x)
&:=   (-1)^m  \frac {(2m + 1)!}{m!} \frac {x}{\mu + \half}
       \Phi(-m,\mu + \frac 32, x^2) \tag2.1.2   \\
&= (-1)^m \frac{(2m+1)!}{m!} \sum _{k=0}^m (-1)^k
       \binom mk \frac {\Gamma(\mu + \half)}{\Gamma(k + \mu + \frac
3 2) }
           x^{2k+1} .        \endalign                      $$
\endproclaim
\proclaim{\jump  2.2 On $ H_{n}^\mu,\,  \gamma_\mu,$  and 
 $ {\bold e}_\mu $}
$$ \align   H_{2m}^\mu(x)
&=    \frac {\Gamma(\mu + \half)}{\Gamma(\half)}
 \frac {\Gamma(m + \half)}{\Gamma(m + \mu +
\half)}{\bbH}_{2m}^\mu(x)
         \tag2.2.1 \\
&= (-1)^m (2m)! \frac {\Gamma(\mu + \half)}{\Gamma(m + \mu +
\half)}
      L_m^{\mu - \half}(x^2) .    \endalign     $$
$$  \align H_{2m+1}^\mu(x)
&= \frac {\Gamma(\mu + \half)}{\Gamma(\half)}
  \frac {\Gamma(m + \frac 3 2)}{\Gamma(m + \mu + \frac 3 2)}
  {\bbH}_{2m+1}^\mu(x)
    \tag2.2.2  \\
&= (-1)^m (2m+1)! \frac {\Gamma(\mu + \half)}{\Gamma(m + \mu +
\frac 3 2)}
     x L_m^{\mu + \half}(x^2)  .    \endalign $$
\endproclaim

    We list the first few generalized Hermite polynomials:
$H_0^\mu(x) = 1$ ,
$$H_1^\mu(x) = (1 +  2\mu )^{-1}2x  ,\qquad H_2^\mu(x) =
  (1 + 2\mu)^{-1}4x^2  - 2       $$
$$H_3^\mu(x) =  (1 + 2\mu )^{-1}(3 + 2\mu )^{-1} {24x^3}
    -(1 + 2\mu )^{-1}12x  $$
$$ H_4^\mu(x) =
\bigl ( (1 + 2\mu)(3 + 2 \mu)\bigr)^{-1 } {48x^4}
 -(1 + 2\mu)^{-1} {48x^2} + 12.                           $$

This class of generalized Hermite polynomials has
a rather nice generating function formula involving the confluent
hypergeometric function $\Phi$ .  If
$\mu \in {\Cal C}_o $ we define
$$\split    {\bold e} _\mu(x)
 & := e^x \Phi(\mu,2\mu + 1, -2x),\qquad \text{ which} \\
 &= \Gamma (\mu + \half) (2/x) ^{\mu - \half}
 \bigl( I_{\mu - \half}(x) + I_{\mu + \half}(x) \bigr)  \\
 &={(2x)^{\half - \mu}}  M_{-\half, \mu} (2x)
                    ,  \endsplit    \tag2.2.3  $$
 where $ I_\mu$ is the modified  Bessel function and $  M_{..} $
is the Whittaker function.
$ {\bold e}_{\mu} $ plays the role of  a generalized exponential
function in what
follows, and indeed $ {\bold e} _0 (x) = e^x $. $ \bold e_\mu $ is
an entire function, say,
$$ {\bold e}_\mu(x) = \sum_{m=0}^{\infty}
\frac{z^m}{\gamma_\mu(m)},
                                             \tag 2.2.4 $$
where  the power series representation for the associated
Bessel function yields
$$ \gamma_\mu(2m) = \frac{2^{2m} {m!} \Gamma(m + \mu +
\half)}{\Gamma(
     \mu + \half)}
   = (2m)! \frac{\Gamma(m + \mu + \half)}{\Gamma(\mu + \half)}
\frac{\Gamma(\half)}{\Gamma(m  + \half)},\, \text { and } \, \tag
2.2.5 $$
$$ \gamma_\mu(2m+1) = \frac{2^{2m+1} {m!} \Gamma(m + \mu + \frac 3
2)}
    {\Gamma( \mu + \half)}
   = (2m+1)! \frac{\Gamma(m + \mu + \frac 3 2)}{\Gamma(\mu +
\half)}
         \frac{\Gamma(\half)}{\Gamma(m + \frac 3 2)}  . \tag 2.2.6
$$
\endgraf
$\gamma_\mu$ plays the role of a generalized factorial.
We list a few of the $\gamma_\mu$ :
$\gamma_\mu(0) = 1, $ $\gamma_\mu(1) =1 +  2\mu ,$
 $\gamma_\mu(2) =(1 +  2\mu)2 , $ $\gamma_\mu(3) =(1 +  2\mu)2(3+
2\mu),  $
 $\gamma_\mu(4) =(1 +  2\mu)2(3+ 2\mu)4  $, and
 $\gamma_\mu(5) =(1 +  2\mu)2(3+ 2\mu)4(5+ 2\mu).  $
We note the recursion relation for the $\gamma_\mu$ :
$$ {\gamma}_\mu(n+1) = (n + 1 + 2\mu \theta_{n+1}){\gamma}_\mu(n),
\qquad  n \in {\Bbb N},                       \tag 2.2.7    $$
where $ {\theta}_{n + 1}$ is defined to be 0 if n + 1 is even and
1 if
n + 1 is odd.
    It follows from (2.1.1) and (2.1.2) that for all $ n \in {\Bbb
N} $
$$  H_{n}^\mu(x)  = n! \sum _{k=0}^{[\tfrac{n}{2}]}
 \frac  { (-1)^k (2x)^{n-2k}}  { k! \gamma_\mu(n-2k)} .  \tag2.2.8 
$$
We note that the coefficient of $x^n$ in the expansion of $ H
_n^\mu$
is $2^n n!/ \gamma_\mu(n)$.  We set down for later
reference  integral expressions for the beta function $B(\cdot
,\cdot ).$

\proclaim{\jump  2.3 Lemma }i) Suppose $\mu > 0 \, , \alpha >
-\half 
, \, n \in \Bbb N , \, \text{ and} \,  \, x \in \Cal C $.
  $$ \gather  B( n +  \alpha + \tfrac 1 2, \mu)   =
 \int_{-1}^1 { t^{2n} |t|^{2\alpha}  ( 1-t)^{\mu - 1} (1+t)^{\mu} 
} dt
                           \tag2.3.1 \\
 B(n + \alpha +\tfrac {3}{2}, \mu)    =
 \int_{-1}^1 { t^{2n + 1} |t|^{2\alpha}(1-t)^{\mu - 1} (1+t)^{\mu} 
} dt
                             \tag2.3.2  \\
  \frac{\gamma_{\alpha}(n)} { \gamma_{\mu+\alpha}(n)} =
\frac 1 {B(\alpha + \half , \mu )} \int_{-1}^1 { t^n |t|^{2\alpha}
  (1-t)^{\mu - 1} (1+t)^{\mu}  } dt                \tag2.3.3 \\
 \bold e_{ \mu + \alpha }(x)      =
\frac {1 }{B(\alpha + \half ,\mu )} \int_{-1}^1
 {\bold e}_{\alpha}(xt) {\vert t \vert }^{2 \alpha}
   (1-t)^{\mu - 1} (1+t)^\mu  dt    \tag 2.3.4 \\
  {\bold e}_\mu (x)        =
\frac {1 }{B(\half ,\mu )} \int_{-1}^1
   e^{xt}  (1-t)^{\mu - 1} (1+t)^\mu \, dt \tag 2.3.5 \\
  \endgather      $$ 
$$ \multline ii) \text{ If now } \, 0 < \mu < \half , \\
   e^x     =  \frac{1}{ B(\half - \mu, \mu) } \int_{-1}^1
{\bold e}_{- \mu} (xt)  {\vert t \vert }^{-2 \mu}
   (1-t)^{\mu - 1} (1+t)^\mu  dt  \endmultline   \tag 2.3.6 $$
 $$ \multline iii) \text{ Suppose }  \,  \alpha > -\half, \,
 \mu > -1  \,  \text{  and  } \,  \mu + \alpha > -\half.\text {
Then}
  \,  \bold e_{ \mu + \alpha }(x)   =  \bold e_{\alpha}(x) \,  + 
\\ 
\frac { \mu}{ \mu + \alpha + \half}  
  \frac {1 }{B(\alpha + \half ,\mu + 1 )} \int_{-1}^1
\bigl( {\bold e}_{\alpha}(xt) - {\bold e}_{\alpha}(x)  \bigr)
 {\vert t \vert }^{2 \alpha}
  (1-t)^{\mu - 1} (1+t)^\mu  dt   \endmultline  \tag 2.3.7  $$
\endproclaim
\demo{Proof} Start with the usual integral representation for
 the beta function, $ B(x,y) \mathbreak
 =  \int_0 ^1 t^{x-1} (1-t)^{y-1} \, dt. $
 Use this to  then derive (2.3.1) and (2.3.2).
(2.3.3) and (2.3.4) follow from this and (2.2.5), (2.2.6), (2.2.4).
Set $ \alpha = - \mu $ in (2.3.4) to get (2.3.6)
(2.3.7) is obtained from (2.3.4) by analytically continuing $\mu$.
We use the functional equation $ (x + y) B(x, y + 1) =
y B(x, y) $ for the beta function.  The rest follow easily.
\enddemo
Next we associate with the generalized exponential function ${\bold
e}_\mu $
a generalized derivative operator  ${\goth D}_\mu$ . These objects
are special cases of functions and operators set down by C.F. Dunkl
in his work on root systems associated with finite reflection
groups.
The papers [11] , [12] are particularly relevant here. \endgraf For
the sake
of simplicity we study the action of ${\goth D}_\mu$ on entire
functions. 
\proclaim{\jump  2.4 Definition } i) The linear operator
 ${\goth D}_\mu$ is defined on all entire
functions  $\phi$ on $\Cal C $ by
$$ {\goth D}_\mu \phi(x) = \phi^\prime(x) +  \frac { \mu}{x}
\bigl (\phi(x) - \phi(-x) \bigr ) , \, x \in \Cal C .  \tag2.4.1  
$$
We use the notation $ {\goth D}_{\mu,x}$ when we wish to
emphasize  that ${\goth D}_\mu$  is acting on
funct- \linebreak
ions  of the  variable  $x$.
Thus $ {\goth D}_{\mu,x}(\phi (x)) := ({\goth D}_{\mu}\phi)(x) $.
\newline
ii) $\frak Q$ is defined on all functions $\phi$  on
 $\Cal C$ by  $$ {\goth Q} \phi(x) = x \phi(x). \tag2.4.2   $$
\endproclaim
\proclaim{\jump  2.5 Properties of ${\goth D}_\mu, \,  {\bold
e}_\mu,$
and $ H^\mu_n $ }
Suppose $ \mu \in {\Cal C}_o, \,  n \in \Bbb N , \mathbreak
 x, z \in \Cal C \,
\text{  and } \,  \phi , \, \psi$ are entire functions.
$$  \gather
 ({\goth D}_\mu^2 \phi)(x) = \phi^{\prime \prime}(x) + \frac
{2\mu}{x}
\phi^{\prime}(x) - \frac{\mu}{x^2} \bigl(\phi(x) - \phi(-x) \bigr)
\tag2.5.1
\\  {\goth D}_\mu^j :  x^n \longmapsto \frac{\gamma_\mu(n)}
    {\gamma_\mu(n-j)} x^{n- j}, j = 0, 1, \dots , n ; \quad
 {\goth D}_\mu^j :  1  \longmapsto 0     
           \tag2.5.2  \\
\intertext {If  $ \psi$ is  an even entire function, then }\\
 {\goth D}_\mu (\phi \psi) = {\goth D}_\mu(\phi)\psi +
 \phi{\goth D}_\mu( \psi)       \tag2.5.3 \\
 H_n^\mu(x) = \frac {1}{B(\half,\mu)} \int_{-1}^{1} H_n(xt)
                (1 - t)^{\mu - 1}(1 + t)^\mu dt
  \qquad  \text{if} \qquad   \mu > 0 .    \tag2.5.4 \\
\intertext { \it Exponential property of $ {\bold e}_\mu$ : } \\
  {\goth D}_{\mu,x} : {\bold e}_\mu(\lambda x) \longmapsto
                       \lambda {\bold e}_\mu(\lambda x)      
\tag2.5.5  \\
\intertext {\it Differential equation for $ {\bold e}_\mu $ :} \\
  x {\bold e}_\mu^{\prime \prime}(x)
+ (1+2\mu){\bold e}_\mu^\prime(x) - (1+x){\bold e}_\mu(x) = 0 
\tag2.5.6 \\
 {\bold e}_{\mu}^{\prime \prime}(x) = {\bold e}_{\mu}(x)
- \frac{ 2 \mu}{ 2 \mu + 1} {\bold e}_{\mu + 1}(x) \tag2.5.7 \\
\intertext {\it Generating function for the $ H_n^\mu$ :} \\
 \exp(-z^2) {\bold e}_\mu (2xz) =  \sum_{n=0}^{\infty} H_n^\mu(x)
         \frac{z^n}{n!} ,\   \  \mu \in {\Cal C}_o    \tag2.5.8
\endgather   $$
\endproclaim
\demo{Proof} (2.5.1) to (2.5.3) follow from the definition in
(2.4.1), and
also (2.2.7).
  To prove (2.5.4) substitute the expansion for $H_n = H_n^0 $
provided
by 2.1  in the right side of (2.5.4) and then employ
(2.3.3) with $\alpha = 0$, so $\gamma_\alpha (n) = n!$. 
\endgraf We prove (2.5.5) using (2.5.2) . $ {\bold e}_\mu (\lambda
x) =
\sum_{j=0}^{\infty} {( \lambda x)^j}/{{\gamma}_\mu (j) }$
is mapped by ${\goth D}_\mu$ to
$\sum _{j=1}^{\infty} { \lambda ^j  x^{j -1} }/{ \gamma_\mu (j-1)
}
= \lambda{ \bold e}_\mu(\lambda x).$ 
\endgraf 
    For (2.5.6) refer to Slater \cite{28,p94}.
\endgraf To prove (2.5.7) first assume $\mu >0$ and check the
result using (2.3.5) and the functional equation for the beta
function. Then analytically continue $\mu$ to $\Cal C_o .$
\endgraf 
    The generating function formula  for the classical Hermite
 polynomials is \linebreak  $ \exp(-z^2 +2xz)  =
 \sum_{n=0}^{\infty} H_n(x) {z^n}/ {n!}. $
Use this result, (2.5.4), and (2.3.7) to prove (2.5.8) for
$\mu > 0$. The result for $ \mu \in {\Cal C}_o $ follows by
analytically continuing $\mu$.
\enddemo
\proclaim{\jump  2.6 Properties of $ H_\mu $ } Suppose
$n \in \Bbb N, \,\,  \lambda , x, z  \in \Cal C, \, |z| < 1. $
$$  \gather
{\goth D}_{\mu , x} :  H_n^\mu(\lambda x) \longrightarrow
 2 \lambda  n H_{n-1}^\mu( \lambda x)      \tag 2.6.1 \\
 {(2{\frak Q} - {\goth D}_\mu)} H_n^\mu =
 \frac{\gamma_\mu(n+1)}{(n+1)\gamma_\mu(n)}
    H_{n+1}^\mu ={ (1 + \frac { 2 \mu \theta_{n+1}}{n+1}})
         H_{n+1}^\mu      \tag2.6.2 \\
\intertext {{\it Three term recursion:} Set $H_{-1}^\mu (x) = 0$.
Then }\\
 2nH_{n-1}^\mu + \frac{\gamma_\mu(n+1)}{(n+1)\gamma_\mu(n) }
      H_{n+1}^\mu = 2{\goth Q} H_n^\mu            \tag2.6.3 \\
{\goth D}_{\mu,x}(e^{-{\lambda}^2 x^2} H_n^\mu(\lambda x)) =
 -\lambda e^{-{\lambda}^2 x^2}
\frac{\gamma_\mu(n+1)}{(n+1)\gamma_\mu(n)}
      H_{n+1}^\mu(\lambda x)                         \tag2.6.4 \\
\intertext {\it Rodrigues formula:}
 (-1)^n e^{\lambda^2 x^2} {\goth D}_{\mu,x}^n e^{-\lambda^2 x^2} =
 \lambda^n  \frac {\gamma_\mu(n)}{n!}
            H_n^\mu(\lambda x) \tag2.6.5 \\
  H_n^\mu
  = \frac{n!}{\gamma_\mu(n)}(2{\goth Q} - {\goth D}_\mu)^n  H_0^\mu
 \tag2.6.6 \\
\intertext {\it Inversion formula:}
  \frac{ (2x)^n }{ \gamma_{\mu} (n) } =
      \sum _{k=0}^{[\tfrac{n}{2}]}
\frac { H_{n-2k}^\mu(x)} {k! (n-2k)!}   \tag2.6.7 \\
\intertext {\it Mehler formula:}
 \sum_{n=0}^{\infty} {\frac{\gamma_\mu(n)}{2^n (n!)^2}}
 H_{n}^\mu(x) H_{n}^\mu(y) z^n  =
\frac {1}{(1-z^2)^{\mu + \half}} \exp \bigl( {-(x^2+y^2)
 \frac {z^2}{1-z^2}} \bigr)
 {\bold e}_\mu(2xy \frac {z}{1-z^2})               \tag2.6.8 \\
 \sum_{n=0}^{\infty} {\frac{\gamma_\mu(2n)}{ (2n)! }
 H_{2n}^\mu(x) \frac { (-1)^n}{ n!} (\frac {z}{2} )^{2n}  =
\frac {1}{(1-z^2)^{\mu + \half}} \exp \bigl( 
 \frac {- x^2  z^2}{1-z^2}} \bigr).
              \tag2.6.9  \endgather $$
\endproclaim
\demo{Proof}  Apply ${\goth D}_{\mu,x}$ to both sides
of (2.5.8) with $'x'$  replaced by $'\lambda x'$ . Then, using
(2.5.5), 
$ \sum_{n=0}^{\infty} {\goth D}_{\mu,x}  H_n^\mu(\lambda x)
{z^n}/{n!}
 = \sum_{n=0}^{\infty}
2 \lambda n H_{n - 1}^\mu(\lambda x) {z^n}/{n!},$ 
and thus (2.6.1) follows. \endgraf
     $ 2x H_{n}^\mu(x) -2n H_{n - 1}^\mu(x)  =
{ (1 +  { 2 \mu \theta_{n+1}}/(n+1)} )
         H_{n+1}^\mu(x) $ follows, upon  separately
 considering the even and odd polynomials, 
 from (2.2.8) and derives (2.6.3).
 From this and (2.6.1) we deduce (2.6.2).
\endgraf
   Use (2.5.3) to infer that ${\goth D}_{\mu,x}(\exp(-x^2)
 H_n^\mu(x))
=\bigl ( (-2x + {\goth D}_{\mu,x})  H_n^\mu(x) \bigr) e^{-x^2}   ,$
  which by (2.6.1) and (2.6.3) equals
the right side of (2.6.4) with $\lambda = 1.$ (2.6.4) for general
$\lambda $ follows easily by the easily derived chain rule
formula ${\goth D}_{\mu,x} : f(\lambda x) \longrightarrow
\lambda ({\goth D}_\mu f)(\lambda x) .$
\endgraf
(2.6.5) is proved by induction using (2.6.1). Note that
 $ H_0^\mu(x) = 1 .$
\endgraf
(2.6.6) is proved by induction using (2.6.2).
\endgraf
 We prove (2.6.7):
 ${\bold e}_\mu (2xz) = e^{z^2} \sum_{n=0}^{\infty}  H_n^\mu(x)
         {z^n}/{n!} $ ,    so 
$ \sum_{n=0}^{\infty} {(2xz)^n}/{\gamma_\mu(n)} \allowmathbreak =
( \sum_{j=0}^{\infty}  {z^{2j}}/{j!})
 \sum_{n=0}^{\infty} H_n^\mu(x)
         {z^n}/{n!}. $ Equate like powers of z to deduce
(2.6.7).
\endgraf
Use (2.2.1), (2.2.2),(2.2.5), (2.2.6) and the bilateral generating
function for the Laguerre polynomials
 Erd\'elyi \cite{15 \rm vol 2, p189} to prove that the left side of
(2.6.8) equals
$$\multline
   \sum_{n=0}^{\infty} \bigl(  {n! \Gamma(\mu + \half)}/
 {\Gamma(n + \mu + \half)} \bigr) {L^{\mu - \half}_n (x^2)}
 {L^{\mu- \half}_n (y^2)} z^{2n}    \\
 +  \sum_{n=0}^{\infty} \bigl(  {n! \Gamma(\mu + \half)}
 /{\Gamma(n + \mu + \tfrac {3}{2})} \bigr) xy  {L^{\mu + \half}_n
(x^2)}
 {L^{\mu + \half}_n (y^2)} z^{2n + 1}   \\
  = \bigl( \frac {\Gamma(\mu + \half) }{(1 - z^2)}  \bigr)
 \exp \bigl(  -(x^2 + y^2) \frac {z^2}{1-z^2 } \bigr)
(xyz)^{-\mu + \half} \biggl( I_{\mu- \half}
\bigl( 2xy \frac {z}{1- z^2 } \bigr)  \\
    +  I_{\mu + \half}(\bigl(2xy \frac {z}{1- z^2 }\bigr) \biggr),
\endmultline
$$
which by (2.2.3) equals the right side of (2.6.8).
\endgraf
 To derive (2.6.9) set $y = 0 $ in Mehler's
formula (2.6.8) and note from Definition 2.1 that
$ H_{2n}^\mu(0) = (-1)^n \frac { (2n)!}{n!} $  and
$ H_{2n +1}^\mu(0) = 0 $ for all $ n \in \Bbb N $.
\endgraf
\enddemo
\proclaim{\jump 2.7 Further properties of $ H_n^\mu $ }  
$$ \exp(-y^2 \goth D_{\mu,x}^2 ) x^n =
 \frac {\gamma_\mu(n)}{n!}
            H_n^\mu(\frac{ x}{ 2y} )  y^n , \quad
0 \ne y \in \Cal C , \quad  n \in \Bbb N .  \tag2.7.1$$
$$ \multline
 \exp (\frac{1}{4} t {\goth D}_{\mu,x}^2 )
\bigl( \exp (-\alpha  x^2) {\bold e}_\mu(2zx) \bigr)
= \\  (1 + \alpha t )^{-\mu - 1/2} \exp(\frac { t z^2 }{1 + \alpha
t})
 \exp ( - \frac { \alpha x^2}{1 + \alpha t} ) {\bold e}_\mu(\frac
{2zx}
{1 + \alpha t} ), \\
\quad  \Re \alpha  > 0, \, \Re t > 0, \,  x \in R, \, z \in \Cal C.
\endmultline \tag2.7.2 $$
\endproclaim
\demo{Proof} We prove (2.7.1). We deduce from (2.5.5) that
$ \exp (-y^2 {\goth D}_{\mu , x}^2 )  {\bold e}_\mu( \lambda x)
= \exp (- {\lambda}^2 y^2 ) {\bold e}_\mu ( \lambda x ) $, so
$$ \exp (-y^2 {\goth D}_{\mu , x}^2 )  \sum _{n=0}^{\infty}
\frac {  { \lambda}^n x^n}{\gamma_\mu(n)}
= \exp (- {\lambda}^2 y^2 ) {\bold e}_\mu ( 2 \lambda y \frac {
x}{2y})
= \sum_{n=0}^{\infty} H_n^\mu(\frac {x}{2y})
\frac {y^n {\lambda}^n }{ n!}.  $$
(2.7.1) follows by equating the coefficients of $ {\lambda}^n $.
\endgraf Assume $ \vert  2 vy  \vert <1 $ and use (2.7.1) and
(2.6.8) to
obtain
$$ \exp(-y^2 \goth D_{\mu,x}^2 ) \exp ( - v^2 x^2 ) {\bold
e}_\mu(2uvx)  =
 \exp(-y^2 \goth D_{\mu,x}^2 )  \sum _{n=0}^{\infty} H_n^\mu(u)
 \frac {(vx)^{n} }{n!} =  $$
$$ \multline    \sum_{n=0}^{\infty} \frac{ \gamma_\mu(n) }{ n! }
H_n^\mu(u)  H_{n}^\mu(\frac {x}{2y}) \frac { y^n v^n }{ n!} =
\sum_{n=0}^{\infty} \frac {\gamma_\mu(n) }{ (n!)^2 2^n} H_n^\mu(u)
H_n^\mu(\frac {x}{2y}) (2yv)^n = \\
\frac {1}{ (1-4 v^2 y^2)^{\mu + \half} }
\exp \bigl(  ( -u^2 - \frac{x^2}{4y^2} )
 \frac {4 v^2 y^2}{1-4v^2 y^2}  \bigr)  {\bold e}_\mu( \frac {
2uvx}
{1 - 4v^2y^2} )   \endmultline $$
 Next set $ t = -4 y^2 ,\, z = uv, \,  \text { and }\alpha = v^2 $.
Thus   (2.7.2) holds at least if $ \vert \alpha t \vert < 1 $.
(2.7.2) follows by analytic continuation for the parameters.
\enddemo
The  identities (2.7.1) and (2.7.2) are generalizations of
classical
Hermite polynomial identities set down, using  Boson calculus
techniques,
 by J.D. Louck \cite{20}.
\endgraf The generalized heat equation problem
$$  \goth D_{\mu,x}^2  \psi (x,t) =
   \frac { \partial  \psi (x,t)}{ \partial t} , \text {  with  }
\psi (x, 0) = \phi (x) , \, t  \ge 0 , \, x \in \goth C , 
\tag2.7.3$$
where $\phi$ is given, has the formal solution
$ \exp(t {\goth D}_{\mu,x}^2 ) \phi (x)  $ . This problem is
related
to the work \cite {6} of Cholewinski and Haimo  and reduces to it
when the function $\phi$ is
assumed to be even. From (2.7.1) we see that if $ n \in \Bbb N$ ,
$$  \exp (t {\goth D}_{\mu , x}^2) x^n =
 \frac {\gamma_\mu(n)}{n!}
            H_n^\mu(\frac{ x}{ 2i \sqrt t })  (i \sqrt t)^n
= \frac {\gamma_\mu(n)}{n!}   \sum _{k=0}^{[\tfrac{n}{2}]}
 \frac  {  x^{n-2k} t^k }  { k! \gamma_\mu(n-2k)}  \tag2.7.4$$
The functions given by (2.7.4) are for even integers n the
generalized  heat polynomials of Cholewinski and Haimo
 \cite{6} and \cite {8,section 14}.\endgraf
 We shall pursue the $ L^2_\mu(R) $ theory for the generalized heat
equation in section 3.
\head 3. The generalized Fourier transform
\endhead
\proclaim{\jump  3.1 Definition } i)  $L_\mu^2(R)$ is  the
Hilbert space of Lebesgue measurable functions  f on R
with $$ {\Vert f \Vert}_{\mu} := \left( \int_{-\infty}^{\infty}
{\vert f(t) \vert }^2 {\vert t \vert}^{2\mu} \, dt  
\right)^ \frac{1}{2}  < \infty. $$

ii)  The \bf  generalized Fourier transform operator \sl
${\Cal F}_\mu,\quad  \mu > - \half$ is defined  on
the linear span $\goth S $ of $ \{ \exp (-x^2) x^n : n \in \Bbb N
\} $
in   $ L_\mu^2(R)$ by
$${\Cal F}_\mu f(x) := \left( 2^{\mu + \half} \Gamma(\mu + \half)
\right )^{-1}  \int_{-\infty}^\infty
{\bold e}_\mu (-ixt) f(t) {\vert t \vert}^{2\mu} dt  \tag3.1.1 $$
\endproclaim
   This transform appears in the physics literature on Bose-like
oscillators,  \cite{24, p294} and \cite{22}, and as a special
case of a general transform in  Dunkl \cite{12}. 
\endgraf  We see from (2.2.3) that if $x$ is real,
$$ {\bold e}_\mu(-ix) =  \Gamma (\mu +\half) 2^{\mu -\half}
    \frac { J_{\mu - \half}(x) - i J_{\mu + \half}(x) }
       { x^{\mu - \half}  } =, {\quad  \roman  say, \quad }
 {\bold c}_\mu(x) - i{\bold s}_\mu(x),
  \tag3.1.2 $$
where $ {\bold c}_\mu  $ is real and even and $ { \bold s}_\mu $
is real and odd. The integral in (3.1.1) is well-defined since 
$$ \gather  \vert {\bold e}_\mu(-ix) \vert \le
 C_\mu \bigl( { \vert x \vert }^{\vert \mu \vert} + 1 \bigr)
\text {  where } C_\mu \in R \, \, \text {if} \,
- \half < \mu < 0,      \tag3.1.3   \\     \text {  and } \,
 \vert { \bold e}_\mu ( -ix) \vert \le 1 \text{ if } \mu \ge 0.
\tag3.1.4  \endgather $$ 
(3.1.3) is proved using the asymptotic formula for Bessel functions
\cite {30, \rm Chapt 7} and (3.1.4) follows
easily from (2.3.5). (2.3.5) implies that $ {\bold e}_\mu(-ix) $ is
a positive definite function of x. Thus the following result
follows. 
In the limiting case $ \mu = 0$ the result degenerates into the
obvious
$ cos^2(x) + sin^2(x) \le  1 $.
\proclaim{\jump 3.2 Remark}  Suppose $ \mu > 0 $. Then
$$ \bigl( J_{ \mu - \half}(x) \bigr)^2 + \bigl( J_{\mu + \half}(x)
\bigr)^2 \le \frac 1 {\Gamma( \mu + \half )^2}
 \vert x/2 \vert^{2 \mu -1 } \tag3.2.1 $$
for  $ x \in R $. Equality holds in (3.2.1) if and only if $ x = 0
$.
\endproclaim
\demo { Proof} 
Use  (3.1.4) and (3.1.2.).
\enddemo
\proclaim{\jump  3.3 Some Fourier transform integrals } 
 Suppose $ \Re \lambda > 0 $.
$$ \int_{0}^{\infty} {\bold c}_\mu(xt)
  \exp(-\half \lambda t^2) t^{2 \mu}\, dt
= \frac{  \Gamma (\mu +\half)}{2  \lambda^{\mu + \half} }
 \exp(-  {x^2}/(4 \lambda))
   \tag3.3.1 $$
$$ \int_{-\infty}^{\infty} {\bold e}_\mu(-ixt) \exp(-\lambda t^2)
|t|^{2\mu} \, dt = \frac{ \Gamma(\mu + \half)}{\lambda^{\mu +
\half}}
\exp(-  {x^2}/(4 \lambda))    \tag 3.3.2 $$
$$ \int_{-\infty}^{\infty} {\bold e}_\mu(-ixt) t^n
 \exp (-\lambda t^2) |t|^{2\mu}\, dt =  (- \frac {i}{2})^n
 \frac{ \Gamma(\mu + \half)}{\lambda^{ \half n + \half + \mu }}
\frac {\gamma_\mu (n)}{n!} \exp(- \frac {x^2} {4 \lambda})
 H_n^\mu(\frac {x}{2 \lambda^\half  })        \tag 3.3.3 $$
$$ \multline \int_{-\infty}^{\infty} {\bold e}_\mu(-ixt)
 {\bold e}_\mu(iyt)  \exp ( -\lambda t^2) |t|^{2\mu}\, dt =
 \frac{ \Gamma(\mu + \half)}{\lambda^{ \mu + \half}}
\exp(- \frac {x^2+ y^2} {4 \lambda})
\bold e_\mu( {xy}/(2 \lambda )), \\
x,y \in \goth C    \endmultline     \tag 3.3.4 $$
$$  \multline
 \int_{-\infty}^{\infty} {\bold e}_\mu(-ixt)
  \exp (-\lambda^2 t^2) H_n^\mu(\beta t) |t|^{2\mu}\, dt \\
 =  (-i)^n  { \Gamma(\mu + \half)}{\lambda^{-2\mu - 1}}
 \left( (\beta/\lambda)^2 - 1 \right)^{n/2}
\exp(-x^2/(4 {\lambda^2}))
  H_n^{\mu}\left(\frac {\beta x}
{2  \lambda (\beta^2 - \lambda^2)^\half }\right), \\
 \roman {if } \quad  {\beta}^2 > {\lambda}^2 > 0.
\endmultline     \tag 3.3.5 $$
$$  \int_{-\infty}^{\infty} {\bold e}_\mu(-ixt)
  \exp ( -\half t^2)  H_n^\mu(t) |t|^{2\mu}\, dt =
 {2^{\mu + \half}}  \Gamma(\mu + \half) (-i)^n
\exp (-\half x^2)  H_n^{\mu}(x)     \tag 3.3.6 $$
\endproclaim
\demo{Proof}(3.3.1) is listed as a Hankel transform in
Erd\'elyi \cite{16 \rm vol 2, p29}. (3.3.2) follows from (3.3.1)
and (3.1.2). Apply $ {\goth D}_{\mu,x}^n $ to both sides of (3.3.2)
and
use the Rodrigues formula (2.6.5) to derive (3.3.3).(3.3.4) follows
from
(3.3.3) by use of (2.5.8) and (2.2.4).
\endgraf
(3.3.5) follows from (3.3.4) by use of the expansions
$ (\exp{y^2}){\bold e}_\mu(i2 \beta y t)  =   \mathbreak
  \sum_{n=0}^{\infty}
    H_n^\mu(\beta  t)(iy)^n/{n!} $ and
$ \exp(-A^2 y^2){\bold e}_\mu(2 B y) = \sum_{n=0}^{\infty}
    H_n^\mu(B/A)(Ay)^n/{n!} $,  \newline  which are a
consequence of (2.5.8).  (3.3.6) follows from (3.3.5) by setting
$\lambda^2 = \half $ and $\beta = 1 $.
\enddemo
\endgraf
 The Hilbert Space $ L_\mu^2(R)$ has the inner product
$ {\langle f,g \rangle}_\mu := \int_{-\infty}^{\infty}
f(t)g^*(t)    {\vert t \vert}^{2\mu} \, dt   $,
where $f,g \in L_\mu^2(R) $
and $g^*$ is the complex conjugate of g. Notice that
 $ {\Vert f \Vert}_{\mu} =  {\langle f,f \rangle}_\mu ^\half ,
f \in  L_\mu^2(R)$.
\proclaim{\jump  3.4 Definition} i) Define the \bf generalized
Hermite functions \sl  $ \phi^{\mu}_n  $ on R   by
 $$  \phi^{\mu}_n(x):=
\biggr( \frac {\gamma_\mu (n)}{ \Gamma (\mu + \half) }\biggr)^\half
 \frac {1} {2^{n/2} n!}
 \exp(-\half x^2) {H}_n^{\mu}(x), n \in {\Bbb N}.\tag3.4.1 $$
\endgraf
ii) Define the operators  ${\bk P}_\mu  , {\bk Q}_\mu , {\bk H}_\mu
  $  and ${\bk A }_\mu $ on the finite span  $\goth S $  of the
 generalized Hermite functions by
 $${\bk P}_\mu \phi(x) :=
 -i ( \phi^{\prime}(x) + \frac {\mu}{x}( \phi(x)
- \phi(-x) )   \text { , }
{\bk Q}_\mu \phi(x) := {\goth Q } \phi(x) = x \phi(x) ,\,  
   \text {  and  }  \tag3.4.2 $$
$$ {\bk A}_\mu := 2^{-\half}({\bk Q}_\mu + i {\bk P}_\mu)  \text {
, }
   {\bk H}_\mu := \half ({{\bk Q}_\mu}^2  + { {\bk P}_\mu }^2 ) 
 . \tag3.4.3 $$
\endproclaim
\proclaim{\jump  3.5  Properties of $\phi_n^\mu$ and $ {\Cal
F}_\mu$}
$$ \{ \phi_n \}_{n \in \Bbb N} \text {  is a complete orthonormal
set
in }  L^2_\mu(R)  .
 \tag3.5.1 $$
$$  {\Cal F}_\mu \phi_n^\mu = (-i)^n \phi_n^\mu , n \in
\Bbb N . \tag3.5.2$$
  {\it   Mehler formula: }
$$ 
\multline
\sum_{n=0}^{\infty}  \phi_{n}^\mu(x) \phi_{n}^\mu(y)
 z^n  =  \\  \frac {1}{\Gamma ( \mu + \half)}
\frac {1}{(1-z^2)^{\mu + \half}} \exp ( {-\half (x^2+y^2)
 \frac {1 + z^2}{1-z^2}})
 {\bold e}_\mu(2xy \frac {z}{1-z^2}) \endmultline       \tag3.5.3$$

\endproclaim
\demo{Proof} Suppose $ u,v \in {\Cal C}$. Then
$$  \multline
\sum_{j,k \in {\Bbb N}} \langle {\phi}_j^\mu,{\phi}_k^\mu
{\rangle}_\mu ( \gamma_\mu (j))^{-1/2} ( \gamma_\mu (k) )^{-1/2}
     {(2^\half  u) ^j ( 2^\half v)^k}  \\
=( \Gamma (\mu + \half) )^{-1}  \int_{-\infty}^{\infty}
 \sum_{j,k \in {\Bbb N}}
{H}_j^\mu (t) {H}_k^\mu (t)  \frac {u^j v^k}{j! k!}
    \exp (-t^2) |t|^{2 \mu} \, dt , \\
 \text{ which by (2.5.8)} \\
= ( \Gamma (\mu + \half) )^{-1}
 \int_{-\infty}^{\infty} \exp(-u^2) {\bold e}_\mu ({2}ut)
\exp(- v^2) {\bold e}_\mu ({2}vt) \exp (-t^2) |t|^{2\mu}
\, dt \\
\text { which by (3.3.4) } \,
=  {\bold e}_\mu(2uv)
= \sum _{j \in {\Bbb N}} {2^j  u^j v^j / {\gamma}_\mu(j) } .
\endmultline $$
Thus  we see that $
{\{ \phi_n^\mu \}}_{n \in {\Bbb N}}$
 is an orthonormal set in $L_\mu^2(R) $ . It is complete by much
the same
argument used to prove that the classical $\mu = 0 $ Hermite 
functions
form a complete orthogonal set in $ L_0^2(R)$, see, for example,
Ahiezer and Glazman, \cite {1}, Chapt 1, paragraph 11.
\endgraf
(2.6.8) implies (3.5.4) and (3.3.6) yields (3.5.3).
\enddemo
  \proclaim{\jump  3.6 Theorem} $  {\Cal F}_\mu $
 is a  unitary transformation on  $  L^2_\mu(R) $
 with eigenvalues  1, -1, i, -i .
 $ {\{ \phi_n^\mu \}}_{n \in {\Bbb N}} $
 is a complete orthonormal set of eigenvectors of $ {\Cal F}_\mu.$
  The inverse Fourier transform is given by :
$${\Cal F}_\mu^* f(x) = \left( 2^{\mu + \half} \Gamma(\mu + \half)
\right )^{-1}  \int_{-\infty}^\infty
{\bold e}_\mu (ixt) f(t) {\vert t \vert}^{2\mu} dt , \quad
  f \in \goth S .
 \tag3.6.1   $$
\endproclaim
\demo{ Proof } The first statement is a direct  consequence of 3.5.
From (3.5.3) we see that $ \Cal F _\mu ^\star \phi_n ^\mu = 
i^n \phi_n^\mu$, and it follows that $(\Cal F_\mu ^\star \psi)(t)
=
(\Cal F_\mu \psi)(-t) $ for any $\psi \in \goth S .$ Then (3.6.1)
follows from (3.1.1)
 \enddemo
 \proclaim{\jump  3.7 More on $\phi_n^\mu$ and $ {\Cal F}_\mu$ }
Define
$ \phi_{-1}^\mu = 0 $ and assume $ n \in \Bbb N . $
$$\align    {\bk A}_\mu \phi_n^\mu   &=
\bigl(  \frac {\gamma_\mu (n)}{\gamma_\mu(n-1)} \bigr)^\half 
\phi_{n-1}^\mu
          \tag3.7.1 \\
 {\bk A}_\mu^\star  \phi_n^\mu             &=
 \bigl(  \frac{\gamma_\mu(n+1)}{\gamma_\mu(n)} \bigr)^\half
 \phi_{n + 1}^\mu
       \tag3.7.2 \\
 {\phi}_n^\mu &= ({\gamma}_\mu (n))^{-1/2 }  {\bk A}^{\star n}
 \phi_0^\mu   \tag3.7.3 \\
\bk H _ \mu \phi_n^\mu =
\half ( {\bk P}_\mu^2 + {\bk Q}_\mu^2 ) \phi^\mu _{n}   &=
\half ( {\bk A}{\bk A}^\star  + {\bk A}^\star {\bk A} ) \phi_n^\mu
=
(n + \mu + \half){\phi^\mu_n} \text{   on }  {\goth S }. \tag3.7.4
\\
\bk J_\mu \phi_n^\mu &= (-1)^n \phi_n^\mu  \tag3.7.5    \\
 i( {\bk P}_\mu {\bk Q}_\mu &- {\bk Q}_\mu {\bk P}_\mu )
 =   ( {\bk A}{\bk A}^\star  - {\bk A}^\star {\bk A} )
 = {\bk I}_\mu + 2\mu {\bk J}_\mu \text{ on } {\goth S}.
\tag3.7.6 \\
 {\Cal F}_\mu^{\star} {\bk Q}_\mu  {\Cal F}_\mu  &=
           { \bk P}_\mu  \text { on } \goth S   \tag3.7.7 
 \endalign $$
   Here ${\bk I}_\mu$ is the identity
operator and  ${\bk J}_\mu$ is the unitary involution defined by
\newline  ${\bk J}_\mu \phi(x) = \phi(-x),\, \phi \in L^2_\mu(R),
\, x \in R.$
\endproclaim
\demo{Proof} By use of (2.5.3) one obtains
$$ {\goth D}_{\mu,x} \bigl( \exp(-\half  x^2) \phi (x)   \bigr) =
 \exp(-\half x^2) \bigl( {\goth D}_\mu  - {\frak Q}_\mu \bigr)
\phi(x) $$
for all smooth $\phi$. Thus (3.7.1) and (3.7.2)  follow from
(2.6.1)
and (2.6.2). (3.7.3) is a consequence of (3.7.2).
\endgraf  $ \bk P_\mu, \, \bk Q_\mu, \, \bk H_\mu $ can be written
in
terms of $\bk A_\mu,  \, \bk A_\mu^\star $, so (3.7.4) and (3.7.6) 
can be derived from (3.7.1) and (3.7.2). (3.7.5) is true since $
\phi_{2n}$
is even and $ \phi_{2n+1} $ is odd.
    Now, the multiplication operator $ {\goth Q}_\mu$  clearly has
a unique closed extension to a selfadjoint operator on the set
$ \{ f \in L_\mu^2(R) : {\goth Q}f \in L_\mu^2(R) \}$, and this
operator we
 name $ {\bk Q}_\mu $. Then $  {\Cal F}_\mu^{\star}  {\bk Q}_\mu 
 {\Cal F}_\mu $ is again a selfadjoint operator. 
One shows that  (3.7.7) is true by  using (3.5.2),
(3.7.1) and (3.7.2).
\enddemo
\endgraf  A class of generalized harmonic and conjugate harmonic
functions and a Hilbert transform operator associated with
the generalized Fourier transform operator ${\Cal F}_\mu $
was introduced and studied by Muckenhoupt and Stein in
\cite{21}. They sketch a proof of the following interesting
generalization of the F. and M. Riesz theorem on the
absolute continuity of analytic measures. If $\mu = 0 $,
then the result is the classical one.
\proclaim{\jump  Theorem 3.8   F. and M. Riesz theorem on absolute
continuity of analytic measures }
Assume that $ \mu \in (-\half , \infty) $ , $ a \in R $ and
 $\nu$ is a complex Borel measure on  $ R_{a}:=  [ a , \infty ) $
such that $ \nu $ is   finite  if $ \mu \in [ 0, \infty ) $
and $ \int_{R_{a}} \bigl( \vert t \vert^{-\mu} + 1 \bigr) | \nu |
(dt)
< \infty $
if $ \mu \in (-\half , 0 ) $ . Assume $\nu $ is an {\it analytic
measure}, that is,
$$ {\int}_{ R_{a} } {\bk e}_\mu (ixt) \nu (dt) = 0  \tag3.8.1  $$
for all real x. Then $ \nu $ is absolutely continuous
with respect to linear Lebesgue measure.
\endproclaim
\demo{Proof} See  \cite{21,p88} for the $ a = 0 $ and $ \mu > 0 $
case. In case $ \mu = 0 $ then the theorem is the classical F. and
M. Riesz theorem, \cite {14,\rm p45}. Our proof consists in
showing that the theorem's hypotheses and  (2.3.5), (2.3.6) imply
that
(3.8.1) holds with $\mu = 0$.
\endgraf Assume first that $ -\half < \mu < 0 $ and (3.8.1) holds.
It follows from (2.3.6) that $\int_{R_a} e^{ixt} \nu(dt) = 0$
for all real x, so $\nu $ is absolutely continuous by the classical
result.
\endgraf Assume next that $ \mu > 0 $. Then an application of
(2.3.7) yields \newline   $ \int_{R_a} \bold e_{\mu_o} (ixt)
\nu(dt) = 0 $
for some $ \mu_o < \mu $, so repeated applications of (2.3.7)
reduce
the problem to the case when $ \mu = 0$.
\enddemo
\endgraf We generalize the Gauss-Weierstrass  operator semi-group
of Hille and Phillips \cite{18,p570} and continue to study the
generalized heat equation (2.7.3) .
\proclaim{\jump  3.9 Definition } For $ t > 0 $  and  $ \mu >
-\half $
the  $ L^2_\mu(R) $  operator ${\bk T}_\mu (t) $  is defined by
$$ \multline    \bigl( {\bk T}_\mu(t) \phi \bigr)(x) :=
 \frac {1}{ (4t)^{ \mu + \half} \Gamma ( \mu + \half ) }
\int_{R} \exp(- \frac {x^2+ y^2} {4t})
\bold e_\mu( \frac{xy}{2t}) \phi (y)
 {\vert y \vert}^{2 \mu} \, dy     , \\
x,y \in R .   \endmultline     \tag 3.9.1 $$
\endproclaim
\proclaim{\jump 3.10 Theorem } $$ {\bk T}_{\mu}(t) \phi =
 \exp (-t {\bk P}_{\mu}^2 ) \phi \text { for every }
t > 0 \text{  and } \phi \in L^{2}_{\mu}(R). \tag3.10.1 $$
\endproclaim
\demo{ Proof } Assume $ t > 0 $ and set $u(\cdot , t) =
 \exp(-t {\bk P}_{\mu}^2)\phi$, which by (3.7.5) equals \newline
$ {\Cal F}_\mu^{*} \exp (-t {\bk Q}_\mu^2)  {\Cal F}_\mu $  .
Then if $ \phi \in L^2_{\mu}(R) $,
$$ \multline u(x,t) = \\
  \bigl( 2^{\mu + \half} \Gamma ( \mu + \half) \bigr)^{-2}
\int_{-\infty}^{\infty} {\bold e}_\mu(ix \tau) \exp (-t {\tau}^2 )
 \vert \tau \vert^{2 \mu}
 \bigl(  \int_{-\infty}^{\infty} {\bold e}_\mu(-i y \tau) \phi(y)
 \vert y \vert^{2 \mu} dy  \bigr) d \tau = \\
  \bigl( 2^{\mu + \half} \Gamma ( \mu + \half) \bigr)^{-2}
\int_{-\infty}^{\infty}  \bigl(  \int_{-\infty}^{\infty}
 {\bold e}_\mu(-i y \tau)  {\bold e}_\mu(ix \tau) \exp (-t {\tau}^2
)
 \vert \tau \vert^{2 \mu} d \tau \bigr) \phi(y) \vert y \vert^{2
\mu} dy \\
\text { which by (3.3.4) }
=  \bigl( {\bk T}_\mu(t) \phi \bigr)(x)
\endmultline $$
\enddemo
\proclaim{\jump  3.11 Theorem }
 $$  \multline     {\bk T}_{\mu}(t/4) :
 \bigl( \exp (-\alpha  x^2) {\bold e}_\mu(2zx) \bigr) 
\longrightarrow
 \\ (1 + \alpha t )^{-\mu - 1/2}  \exp(\frac { t z^2 }{1 + \alpha
t})
 \exp ( - \frac { \alpha x^2}{1 + \alpha t} ) {\bold e}_\mu(\frac
{2zx}
{1 + \alpha t} ), \\
 \, \Re t > 0, \, z  \in \Cal C,\,  \Re \alpha > 0 .
 \endmultline \tag3.11.1 $$
$$ \multline
     {\bk T}_{\mu}(t) :
  \exp (-\alpha  x^2)   \longrightarrow  (1 + 4 \alpha t )^{-\mu -
1/2} 
 \exp ( - \frac { \alpha x^2}{1 + 4 \alpha t} ) , \\
 \, \Re t > 0, \, \Re \alpha > 0 .
\endmultline  \tag3.11.2
   $$
\endproclaim
\demo{ Proof}(3.11.1) is a consequence of (2.7.2) and (3.11.2) is
obtained by setting $z = 0 $.
\enddemo
\proclaim{\jump 3.12 Theorem } Suppose 
$ \phi \in L^2_\mu(R) $. Then $ \psi(x,t) := ({\bk T}_\mu
(t)\phi)(x) $
satisfies $$
 \pd{\psi}{t} = - {\bk T}_{\mu}^2 \psi =
  { \goth D}_{\mu }^2 \psi  \in L^2_{\mu}(R)  \tag 3.12.1  $$
$$ \text {  and }  \lim_{t \searrow 0}
 \Vert \psi(\cdot , t) - \phi \Vert_{ \mu} = 0.  \tag3.12.2
 $$
\endproclaim
\demo{ Proof}  By (3.7.5) $ \psi (\cdot ,t) =
 \exp(- t { \bk P }_{\mu}^2) \phi = {\Cal F}_\mu^{*}
\exp (- t {\bk Q}_{\mu}^2 ){\Cal F}_\mu \phi $.
Since  \newline $ {\bk Q}_{\mu}^2 \exp( - t {\bk Q}_\mu^2 ) \phi
\in L^2_{\mu}(R) $ for every $\phi \in L^2_{\mu}(R) $, it follows
that $  {\bk P}_{\mu}^2 \exp( - t {\bk P}_\mu^2 ) \phi \in
L^2_{\mu}(R) $ for every $\phi \in L^2_{\mu}(R) $.
Thus we see that (3.12.1) holds. 
$$ \multline
 \Vert \psi(\cdot ,t) - \phi \Vert_{ \mu}^2 =
\Vert {\Cal F}^\star \exp ( - t {\bk Q}_\mu^2 ) \Cal F \phi
- \phi \Vert_{ \mu}^2 = \Vert \exp( -t {\bk Q}_{\mu}^2 )
\Cal F \phi - \Cal F \phi \Vert_{ \mu}^2  \\
= \int_{- \infty}^{ \infty} \vert \exp ( -t  y^2 ) -1 \vert^2
\vert (\Cal F \phi)(y) \vert^2 \vert y \vert^{2 \mu}  dy
 \longrightarrow 0 \text{ as } y \searrow 0 .\endmultline  $$
 (3.12.2) follows.
\enddemo
\proclaim{\jump  3.13 Corollary} Assume the hypotheses and notation 
of
Theorem 3.12 , so $ \psi(t, x) $ satisfies (3.12.1) and (3.12.2).
i) Suppose that $ \phi $ is also an even function. Then
$$  \pd{\psi}{t} = \dfrac {\partial^2 \psi}{\partial x^2}
+ \frac {2 \mu}{x}\pd{\psi}{x}                   \tag3.13.1 $$
ii) Suppose that $ \phi $ is also an odd function. Then
$$    \pd{\psi}{t} =
 \dfrac {\partial^2 \psi }{\partial x^2}
+ \frac {2 \mu}{x} \pd{\psi}{x} - \frac{2 \mu}{x^2} \psi
 \tag3.13.2  $$
\endproclaim
\demo{ Proof } Use (2.5.1) and (3.12).
\enddemo
\head 4. Generalized translation
\endhead
\proclaim{\jump  4.1 Definition } i) The {\bf generalized
translation
 operator} \newline
 ${\goth T}_y,\, y \in R $  is defined  by
$$  {\goth T}_y \phi  := {\bold e}_\mu(y {\goth D}_\mu) \phi
  = \sum_{n=0}^{\infty} \frac {y^n}{\gamma_\mu(n)}
{\goth D}_\mu^n \phi
$$ for all entire  functions $\phi$ on $\Cal C $ for which the
series converges pointwise.
\endgraf
ii) The linear operator ${\Cal T}_y,\,  y \in R $ is defined on
$L_\mu^2(R)$ by
$$ { \Cal T}_y \phi := {\bold e}_\mu(iy {\bk P}_\mu) \phi .
            \tag4.1.1  $$
We use the notation ${\goth T}_{y,x}$ 
when we wish to emphasize the
functional dependence on the variable $x$. \newline
iii) The     $\mu -$\bf binomial coefficients \sl are defined by
$$ {\binom nk}_\mu := {\gamma_\mu(n)}/\bigl(\gamma_\mu(j)
 \gamma_\mu(n-j)\bigr), \quad  k = 0, \dots , n, \, n \in \Bbb N.
\tag4.1.2 $$
iv) The $\mu -$\bf binomial polynomials
 \sl $ \{ p_{n , \mu }(\cdot,\cdot) \}_{n \in \Bbb N } $
are defined by
$$ p_{n , \mu}(x,y) := {\goth T}_{y,x} x^n 
  = \sum_{j=0}^{\infty} \frac {y^j}{\gamma_\mu(j)}
{\goth D}_{\mu,x}^{j} x^n .  \tag4.1.3 $$
\endproclaim
Notice that $ {\Cal T}_y = {\goth T}_y $  for almost all real $y$
on the class  of $L^{2}_\mu (R) $ entire functions of the
form 
$\{ p(x) \exp(- \lambda x^2) : p \, \, \text {is a complex
polynomial and}
\,   \lambda > 0 \}.$  
\proclaim{\jump  4.2   Properties of ${\goth T}_y$    }
  { \it  $\mu$-binomial expansion :}
$$   p_{n , \mu}(x,y) =   \sum_{j=0}^{n}
 {\binom nk}_\mu  \, x^j y^{n-j} .          \tag4.2.1 $$
 The first few $\mu-$binomial polynomials are
$ 1, \quad x + y, \quad  x^2 + \frac 2 {1 +2 \mu} xy
+ y^2 $ , and $ \quad x^3 + \frac{ 3 + 2 \mu}{ 1 + 2 \mu} ({x^2}y
+ x{y^2}) + y^3 , \quad x^4 + 4 \frac {1}{ 1 + 2 \mu } ({x^3}y +
x{y^3})
+ 2 \frac{ 3 + 2 \mu}{1 + 2 \mu}x^2 y^2 + y^4 .  $
Thus it is clear that  ${\goth T}_y $ in general does not take
nonnegative functions into nonnegative functions. Consider, for
example,
$ (x - 1)^2 .$ 
$$   {\goth T}_{y,x}{\bold e}_\mu(\lambda x) =
{\bold e}_\mu(\lambda y) {\bold e}_\mu(\lambda x), \quad
\lambda \in {\Cal C} \tag4.2.2 $$
{\it Generating function : }
$$  {\bold e}_\mu( \lambda x )  {\bold e}_\mu( \lambda y ) =
\sum_{n=0}^{\infty} \frac { p_{n , \mu}(x,y)}{ \gamma_\mu(n)}
 {\lambda}^n , \quad
   \text{ if } \quad \lambda \in \Cal C .  
\tag4.2.3 $$
$$ \align   {\Cal T}_y : \exp (- \lambda x^2) \longmapsto &
 \exp(- \lambda ( x^2 + y^2))
      {\bold e}_\mu (-2 \lambda xy)
\text {\quad if \quad}    \lambda > 0   
 \tag4.2.4 \\
   {\Cal T}_y : x \exp (- \lambda x^2) \longmapsto &
(x + y) \exp(- \lambda (x^2+ y^2))
     {\bold e}_\mu (-2 \lambda xy) 
\text { \quad if \quad }  \lambda > 0
 \tag4.2.5\endalign $$
\endproclaim  
\demo{ Proof } Use (2.5.2) so
$$ p_{n , \mu}(x,y) =
 \sum_{j=0}^{\infty} \frac{ y^j }{ \gamma_\mu (j)}
\goth D_{\mu,x}^j  x^j  = \sum_{j=0}^{n}
 {\gamma_\mu(n)}/\bigl(\gamma_\mu(j)
 \gamma_\mu(n-j)\bigr) x^j y^{n-j}
,$$  which is equivalent to (4.2.1).
 $$\multline    {\goth T}_{y,x} {\bold e}_\mu(\lambda x)  
 = {\bold e}_\mu(y \goth D_{\mu,x}){\bold e}_\mu(\lambda x) 
  = \\    \sum_{n=0}^{\infty} \frac {y^n}{\gamma_\mu(n)}
{\goth D}_{\mu,x}^{n} {\bold e}_\mu (\lambda x) = 
 \sum_{n=0}^{\infty} \frac {\lambda^n   y^n}{ \gamma_\mu(n) }
\bold e_\mu(\lambda x)
= \bold e_\mu( \lambda y) \bold e_\mu(\lambda x) \endmultline $$ 
This proves (4.2.2). \newline 
\quad  (4.2.3) is a consequence of (4.2.2) and (4.1.3) .  
(4.2.4) follows from (2.6.5),(2.5.8) since:
$$ \multline  \text { Set } \lambda = \alpha^2. \text { Then }
 {\goth T}_y : \exp (- \alpha^2 x^2) \longmapsto
 {\bold e}_\mu(y {\goth D}_{\mu,x})\exp(- \alpha^2 x^2) \\
=  \exp ( - \alpha^2 x^2 ) \sum_{n=0}^{\infty}  \frac
{y^n}{\gamma_\mu(n)}
 { \goth D}_{\mu,x}^n   \exp (- \alpha^2 x^2)
 = \sum_{n=0}^{\infty} (-1)^n \frac{ \alpha^n y^n }{n!}
 H_n^\mu(\alpha  x) \\  = \exp ( - \alpha^2 x^2 )
\exp ( - \alpha^2 y^2) {\bold e}_\mu(-2 \alpha^2 xy )
\endmultline $$
\quad (4.2.5) follows when one applies $ {\goth D}_{\mu, x } $
to both sides of (4.2.4).
\enddemo
In the rest of section 4 we assume $ \mu > 0 $.
 Define the probability  measure $ \alpha_\mu $ and the function
 $ \widetilde \omega $    by
$$ \alpha_\mu(dt) := \frac {1}{ B (\half, \mu) }
 { (1  - t)^{\mu - 1}(1 + t )^{\mu} }
 \, dt,    \quad  t \in ( -1, 1 ), \tag4.2.6 $$
 $$  \widetilde \omega (t) := ( x^2 + 2xyt + y^2 )^{\half},
 \quad t \in [-1,1]  \text{  and }  x,y \in R .  \tag4.2.7 $$

\proclaim {\jump  4.3 Lemma   }
Suppose   $ \phi $ is an $ L^{\infty}(R) $ function. Then  
 $$ \multline {\goth T}_{y , x}( \phi(x) ) =
\frac {1}{ 2} {\int}_{-1}^{1} \bigl( 1 +
\frac { x + y }{ \widetilde \omega (t) } \bigr)
\phi( \widetilde \omega(t) ) \, \alpha_\mu( dt ) + \\
\frac {1}{ 2 } {\int}_{-1}^{1} \bigl( 1 -
\frac { x + y }{ \widetilde \omega (t)} \bigr) 
\phi(- \widetilde \omega(t) ) \, \alpha_\mu( dt )
\endmultline \tag4.3.1 $$
\endproclaim
\demo{ Proof } It follows from (4.2.4),(4.2.5) and (2.3.5) that 
$$  (\goth T_y \phi)(x) = \int_{-1}^{1}
\psi \bigl( \widetilde \omega (t) \bigr)  \, \alpha_\mu(dt), \qquad
  (\goth T_y \psi)(x) =( x + y)\int_{-1}^{1}
 \frac {1}{\widetilde \omega (t)}
  \psi \bigl( \widetilde \omega (t)\bigr) \, \alpha_\mu(dt) $$
 if $ \phi $  is an even and $\psi$ is an odd $ L^{\infty}(R)$ 
function. 
 Thus (4.3.1) follows.
\enddemo
\proclaim{\jump 4.4 Corollary } Suppose $ \mu >  0 $
and $\phi$ is an  entire function. Then
$$ \align  ({\goth T}_y \phi)(x) &= ({\goth T}_x \phi)(y) \text{
for  all
real x and y }. \\
  {\goth T}_{y,x}\bigl( x  \phi(x) \bigr) &= (x + y)
  {\goth T}_{y,x}( \phi(x) ) \quad \text {if $ \phi$ is even }.\\
 p_{2n + 1, \mu}(x,y) &= (x + y) p_{2n, \mu}(x, y) ,   \, n \in
\Bbb N .
\endalign $$
\endproclaim
\demo{Proof} These follow from (4.3.1).
\enddemo
\proclaim{\jump 4.5 Notation } i) Suppose $ x , y,  \xi \in R  $.
Define
$$ \align
 \Psi(x,y, \xi ) &:= \frac{1}{16} \bigl( (x + y)^2 - {\xi}^2)
 ({\xi}^2 - (x -y)^2 \bigr)   \\
 \Xi &:= \{ (x,y,\xi) \in R^3 : \Psi (x, y, \xi) > 0  \}
\tag4.5.1\\   
 \Xi (x,y) &:= \{ \xi \in R : \Psi ( x, y, \xi) > 0 \}
  \endalign $$
Thus
$$ \Xi(x , y) = \bigl(- \vert x + y \vert , - \vert x - y \vert
\bigr) \cup \bigl( \vert x - y \vert , \vert x + y  \vert \bigr)
 \text { if } \, xy > 0  $$
$$ \Xi( x , y) = \bigl( - \vert x - y \vert ) , - \vert x + y \vert
\bigr) \cup \bigl( \vert x + y \vert , \vert x -y \vert \bigr)
 \text{ if } \, xy < 0 . $$
 ii) We note that  $\Psi$ is the homogeneous symmetric polynomial
 that is relevant in Heron's formula for the area of a triangle,
and
$\Psi(x,y,\xi) = \Psi(\vert x \vert, \vert y \vert, \vert \xi
\vert). $
It appears in the expressions for the generalized translation
operator
that form a basis for the Bessel calculus studied  by 
Cholewinski \cite{8} and by  others. See \cite{2,p35-36}
and \cite{8}  for references but note that  operator they study
 acts only on  even  functions or on functions on a half-line.
The translation operator we will set
down acts  more  generally on functions on $R$.     
 iii) Suppose $ x , y , \xi \in R $.
Define $ \Delta (x , y, \xi ) $ to be the area of the triangle
formed, if possible, with sides of length $ \vert x \vert ,
\vert y \vert, \vert \xi \vert , $ \, and
0 otherwise. Then
Heron's formula, see \cite {9,p12} states that
$$ \align  \Delta(x ,y, \xi) &=  (\Psi(x,y, \xi )^\half \quad
  \text { if } \, (x, y, \xi) \in \Xi  \text{ and } \tag4.5.2 \\
\Delta(x , y, \xi) &= 0  \quad  \text{ if } 
 (x , y,\xi) \notin \Xi 
\endalign     $$ 
iv) Define the measure $ \beta_{\mu,x,y} $  on $\Xi(x,y) $   by
$$  \beta_{\mu, x ,y} (d \xi) := \frac {1} { B(\half,\mu) }
\bigl(  \frac { \Delta(x, y, \xi) }{ \vert xy \vert} \bigr)^{2 \mu}
  \, d \xi .
 \tag4.5.3  $$ \endproclaim
\proclaim{\jump 4.6 Lemma } Assume $ x ,y \in R \setminus \{ 0 \}
$
,  $\phi, \psi \in L^{\infty} (R) $ with $\phi$ an even and
$\psi$  an odd function. Then
$$  \int_{ \Xi (x , y) } \frac{ sgn(xy \xi)}{ x + y - \xi}
  \phi ( \xi)  \,  \beta_{\mu, x, y }( d \xi).
 =  \int_{-1}^{1} \phi(\widetilde \omega (t) ) \alpha_\mu (dt) 
 \tag4.6.1  $$
$$\int_{ \Xi (x , y) } \frac{ sgn(xy \xi)}{ x + y - \xi}
  \psi ( \xi)  \,  \beta_{\mu, x, y }( d \xi)=
(x + y) \int_{-1}^{1} \frac{  \psi(\widetilde \omega (t) )}
{\widetilde \omega (t)} \alpha_\mu (dt) \tag4.6.2  $$.
\endproclaim
\demo{ Proof }  Set $\xi = \tild
 = (x^2 + y^2 + 2xyt)^{\half} , \text{ so }
t = \frac{ {\xi}^2 - x^2 -y^2 }{2xy} $ Thus
$$ 1 + t = ( \xi^2 - (x - y)^2)/ {(2xy)}, \quad
   1 - t = ( (x + y)^2 - \xi^2)/{(2xy)}, \text{ and}  $$
 $$  (1 - t^2)^\mu =
 \bigl(  \frac{2 \Delta(x,y,\xi)}{\vert xy \vert} \bigr)^{2 \mu}.
$$
    First assume $xy > 0$.
 Then $ B(\half, \mu)$ times  the right side of (4.6.1)
$$  = \bigl(  \int_{ -\vert {x+y} \vert } ^{ -\vert {x - y }\vert
} +
  \int_{ \vert {x-y} \vert } ^{ \vert {x + y }\vert }  \bigr)
\frac {sgn (xy \xi)}{x + y - \xi}
\phi (\xi)      \, \beta_{\mu, x , y}( d \xi ) \quad  = $$
$$  \int_{ \vert {x-y} \vert } ^{ \vert {x + y }\vert } 
\frac{ 2 \xi}{ (x+y)^2- \xi^2}
 \bigl( \frac{2 \Delta(x,y,\xi)}{\vert xy \vert} \bigr)^{2 \mu}
   \phi (\xi) d \xi       
 = \int_{-1}^{1} \phi(\tild)  (1 -t)^{\mu-1} (1 + t) ^{\mu}  \, dt. 
$$
This implies (4.6.1) in case $xy > 0$. The case $xy < 0$ follows
 similarly.  
    From (4.6.1) we deduce that the right term in (4.6.2)
equals
$$(x + y) \int_{ \Xi (x , y) } \frac{ sgn(xy \xi)}{( x + y -
\xi)\xi }
  \psi ( \xi)  \,  \beta_{\mu, x, y }( d \xi) $$
$$  =  
\int_{ \Xi (x , y) } { sgn(xy \xi)}\bigl( \frac {1}{ x + y - \xi}
+
 \frac {1}{\xi} \bigr)
  \psi ( \xi)  \,  \beta_{\mu, x, y }( d \xi)= \text{ the left side
  of (4.6.2)}.
$$
\enddemo
\proclaim{\jump  4.7 Theorem } Suppose  $\phi \in L^{\infty}(R)$,
$ \mu > 0 $, and $ x, y \in R \setminus \{0 \}.$
$$  ({\goth T}_y  \phi)(x)  =
 = \int_{ \Xi (x , y) } \frac{ sgn(xy \xi)}{ x + y - \xi}
  \phi ( \xi)  \, \beta_{\mu,x,y}(d \xi)  \tag4.7.1$$
\endproclaim
\demo { Proof} (4.7.1) holds for even and odd  $\phi$ by Lemma 4.6
and
Theorem 4.3.  Thus it holds for all $\phi$. 
\enddemo
\proclaim{\jump  4.8 Theorem } Suppose $ \phi \in L^{\infty}(R) $
and
 $ \mu >  0 . $
$$  ({\goth T}_y \phi)(x) = ({\goth T}_x \phi)(y) \text{ for almost
all
real x and y } \tag4.8.1 $$
$$ \Vert \goth T _y \phi \Vert_{ \mu}  \le
\Vert \phi \Vert_{ \mu} , \phi \in L^{2}_{\mu}(R), y \in R .
\tag4.8.2 $$
\endproclaim
\demo { Proof} (4.8.1) is clearly implied by (4.7.1), and (4.8.2)
by (3.1.4) and (4.1.1).
\enddemo

\head 5. The Bose-like oscillator               
\endhead                             

 Suppose that $\frak H$ is a  complex Hilbert space.
We will be examining certain equations of motion and commutation
relations that relate several unbounded operators. In order
to avoid the pitfalls associated with formal computation
involving unbounded operators \cite{26,\rm p270-274} we shall
postulate the existence of a suitably tailored
linear invariant set of analytic vectors.
\proclaim{\jump Definition 5.1}
i)Suppose  $ \bk P , \bk Q $ and $ \bk H $ are possibly unbounded
selfadjoint operators on $\goth H $ . $ \bk P $ and $  \bk Q $ are
\it dominated \rm by $ \bk H $ if 
$$ \goth S  \text{ is a linear  invariant set of analytic vectors
for }
 \bk P, \bk Q, \text {  and  }  \bk H ,      $$
 and 
 $\goth S :=  \bigcup \{ \bk E ((-n,n)) \goth H :n \in \Bbb N \}$
,
where $ \bk E $ is the spectral measure of $ \bk H $ .
\endgraf
ii) We next specify $\bk H $.
Let  $ V(\cdot)$ be a continuously differentiable real function on
the real
line with derivative $ V^{\prime}(\cdot)$ and specify the
associated
{\it Hamiltonian 
operator }   $\bk H $ by $ \bk H =  \half { \bk P}^2 + V(\bk Q). $ 

Then the  {\it  equations of motion associated with the Hamiltonian
} are 
 $$ i[\bk P, \bk H ] : = i( \bk P \bk H - \bk H \bk P) = 
V^{\prime}(\bk Q)
\text{ \quad and \quad } i[\bk Q, \bk H ] : = i(\bk Q \bk H - \bk
H \bk Q)
 = - \bk P    \tag{$ \Bbb{ EM}$  }$$    on $ \goth S $   .
 In case $ V(\bk Q) = \half {\bk Q}^2 $
the equations of motion are that of the { \it quantum mechanical
harmonic
oscillator. } 
\endgraf
iii) Suppose $ \bk P \, \text{ and } \,  \bk Q $ are selfadjoint
operators
that are dominated by the selfadjoint operator $ \bk H $ , where
$ \bk H := \half (\bk P ^2  + \bk Q ^2 ) $  on $ \goth S $ ,
where $ \goth S $ is as in 5.1i).
Then $ (\frak H, \bk P, \bk Q, \bk H)$ is a  \bf  Bose-like
(quantum mechanical simple harmonic) oscillator   \sl ,
or a \bf para-Bose oscillator \sl if  the equations of motion are
 $$  i [\bk P, \bk H ]  =  \bk Q  \text{ \quad and \quad }
 i [\bk Q, \bk H ] = - \bk P .  \tag 5.1.1  $$ 
\endgraf
iv) The Bose-like oscillator is \bf irreducible  \sl if whenever
 $\bk B$
is an everywhere defined bounded operator on $\frak H $ to
$ \frak H $ such that
$$ \exp(i \lambda \bk P) \bk B = \bk B  \exp(i \lambda \bk P)
\text {\quad and \quad }
\exp(i \lambda \bk Q) \bk B = \bk B  \exp(i \lambda \bk Q)
\tag5.1.2  $$
 for all real  $\lambda$  , then there exists $ c \in \Cal C $
such that $ \bk B = c \bk I $.
 \endproclaim
\endgraf
\proclaim{\jump Remark 5.2} $\bk Q$ is the \bf position \sl, $\bk
P$ is the
\bf momentum \sl, $\bk A :=  2^{-\half} (\bk Q + i \bk P) $ is the
\bf annihilation \sl and ${\bk A}^\star = 2^{-\half} (\bk Q - i \bk
P) $
is the \bf creation \sl operator. $\bk H = \half( {\bk A}^\star \bk
A
+ \bk A {\bk A}^\star ) $.
\endproclaim
\proclaim{\jump Remark 5.3}
  i)By a \it quantum mechanical system consisting of a single
particle
moving in one dimension \rm we mean a triple of self adjoint
operators 
$ \bk P $ and $  \bk Q $ dominated by $ \bk H = \half \bk P^2 +
V(\bk Q)$
  that satisfy the  commutation relation 
$$ i[\bk P, \bk Q] := i(\bk P \bk Q - \bk Q \bk P) =  \bk I
\text {  on   } \goth S \tag{$\Bbb {CR}$} $$
If $\Bbb{CR}$ 
holds, then the equations of motion $ \Bbb {EM} $  holds on $ \goth
S $.
\endproclaim
\demo{Proof}The following computations are valid on $\goth S$ :
$$ i [\bk P,\bk H] =i [\bk P, V(\bk Q)] =  V^{\prime}( \bk Q ),
 \text {  and  }$$  
$$ 2i[\bk Q,\bk H] = i [\bk Q,{\bk P}^2] = i [\bk Q,\bk P] \bk P
+ i \bk P [\bk Q,\bk P] = -2 \bk P, $$
which imply $\Bbb{ EM}$ .  
\enddemo
   In 1950 E. P. Wigner \cite{31} posed the question
\lq Do the equations of motion determine the quantum mechanical
commutation relations? \rq  He considered formal operator
equations of motion of  the form $\Bbb {EM}.$ In the  case where 
the equations of motion are that of a Bose-like oscillator,
Wigner noted that there exists
a one parameter family of inequivalent operator representations
for position $\bk Q$ and momentum $\bk P$. There is now an
extensive
physical literature on these representations, see \cite{24} and
 \cite{22}. We shall give a self-contained treatment
of the theory , and indicate their relationship  to
the  generalized Fourier transform and generalized Hermite
functions.
\endgraf  Our goal is to  study a  generalization of the Boson
calculus.
The  \it {Boson calculus} \rm is the collection of operators,
 functions, and analysis
associated with the quantum mechanical harmonic oscillator. It is
often studied using Lie group theory, by noting that $\Bbb {CR} $
gives
rise to the Heisenberg group through the Weyl commutation
relations.
An alternate perspective to the Boson calculus appears in Glimm and
Jaffe \cite {17,\rm Chapt1} and Biedenharn and Louck
 \cite{3,Chapt5}.
Our perspective uses operator equations and
operator theory directly to implement the analysis, without 
emphasizing the Lie group aspects of the algebraic structure
of the operator equations. \endgraf  The generalization,
the  \it {Bose-like oscillator calculus, }\rm is set down in
Ohnucki and Kamefuchi \cite{24,Chapt23} and in \cite{22}.
For a related calculus see Cholewinski \cite{7}.\endgraf    
        From now on we consider a Bose-like oscillator
and fix the notation of definition 5.1. We do not assume
that $\Bbb{ CR }$ holds. 
\proclaim{\jump  5.4  Alternate formulations of the equations of
motion }
\endgraf  The following statements are equivalent to (5.1.1):
 $$ \alignat 2
  i[\bk P , {\bk Q}^2] &= 2 \bk Q \quad  \text{  and  } \quad
  i[  {\bk P}^2,  \bk Q ] &&= 2 \bk P \text{ on }  \goth S ;
\tag5.4.1 \\
 [ \bk A , \bk H ] &= \bk A \quad   \text{  or  } \quad
                  [ {\bk A}^\star , \bk H ]
    &&= - {\bk A}^\star  \text{ on }  \goth S  ; \tag5.4.2 \\
 [ \bk A , {\bk A}^{\star 2} ]   &= 2 {\bk A}^\star \,  \text{ or 
}
\quad   [{\bk A}^\star , {\bk A }^2 ] &&= - 2 {\bk A}
 \text{ on }  \goth S . \tag5.4.3 \endalignat  $$
\endproclaim
\demo{Proof} (5.4.1)-(5.4.3) follow from (5.1.1) and 5.1, 5.2.
\enddemo
\proclaim{\jump  5.5 Remark}  Suppose $ n \in \Bbb{ N }.$
When considered on $ \goth S $ :
$$ \gather  [\bk A , \bk A ^\star ] = i[ \bk P , \bk Q] 
  \text{ commutes with  }  {\bk A}^2 , {\bk A}^{\star 2},  
 \text {  and  } \bk P^2 , \bk Q^2 ,  \bk H ;  \tag5.5.1  \\ 
 [ \bk A , {\bk A}^{\star (2n) } ] = 2n{ \bk A}^{\star (2n-1)} ,
\tag 5.5.2 \\
 [ \bk A ,  {\bk A}^{\star (2n+1) } ] =
 {\bk A }^{\star (2n) }
( 2n + [ \bk A , \bk A ^ \star ] ) ;\tag 5.5.3  \\
  i[\bk P ,   {\bk Q}^{2n}]  = 2n{\bk Q}^{2n-1}
\text{  , \qquad }
  i[\bk P , {\bk Q}^{2n+1}] = {\bk Q}^{2n} (2n + i [\bk P, \bk Q ])
;
\tag 5.5.4         \\
  i[{\bk P}^{2n} , \bk Q]  = 2n {\bk P}^{2n - 1}
\text{  , \qquad}   
  i[{\bk P}^{2n + 1}, \bk Q] = {\bk P}^{2n} (2n + i [\bk P, \bk Q])
\tag 5.5.5 
 \endgather $$
\endproclaim
\demo{ Proof} (5.5.1) follows from 5.4 and the Jacobi identity
$$ \bigl[ [\bk X, \bk Y], \bk Z \bigr] +
   \bigl[ [\bk Y, \bk Z], \bk X \bigr]
+  \bigl[ [\bk Z, \bk X], \bk Y \bigr] = 0 , $$ which holds
for suitably defined operators $ \bk X , \bk Y, \bk Z $  selected
from $ {\bk A}^2 , {\bk A}^{\star 2}, \bk H  \mathbreak \text 
{and \,}  \bk A , \bk A^\star . $
\endgraf   (5.5.2) is true if $n = 1$ by 5.4.
Assume that it is true for n. Then $ (2n+2){\bk A}^{\star (2n+1)}
=
(\bk A {\bk A}^{\star 2n} -  {\bk A}^{\star 2n} \bk A ) {\bk A }^{
\star 2}
+ 2 {\bk A}^{\star (2n+1)}$ which equals $   \bk A {\bk A}^{\star
(2n+ 2) } -
 {\bk A}^{\star 2n}( \bk A {\bk A}^{\star 2} - 2 {\bk A}^\star ).$
By (5.4.3) this $= [\bk A, {\bk A}^{\star (2n+2}] $. This proves
 (5.5.2).
\endgraf (5.5.3) is clearly true if for $n=0 . $ Assume that it is
true for n. Then
$ [ \bk A , {\bk A}^{\star ( 2n +3 )} ] =
\bigl( \bk A {\bk A}^{\star (2n+1)} - {\bk A}^{\star (2n + 1)} \bk
A
\bigr) {\bk A}^{\star 2 }    + {\bk A}^{\star (2n + 1) }
\bigl( \bk A {\bk A}^{\star 2} - {\bk A}^{\star 2} \bk A  \bigr) .
$
But this equals
 $ {\bk A}^{\star ( 2n + 2 )} \bigl( 2n + [ \bk A , {\bk
A}^{\star}]
 \bigr) +   {\bk A}^{\star (2n + 1)} \bigl( 2 { \bk A}^{\star } 
\bigr), $
so (5.5.3) follows.
\endgraf (5.5.4) and (5.5.5) are proved similarly.
\enddemo
\proclaim{\jump  5.6 Lemma} The following operator identities
hold for all real $\lambda ,\, \nu $.
$$ \gather \exp(i \lambda \bk H) \bk Q \exp(- i \lambda \bk H) =
 \bk Q  cos{ \lambda}  +   \bk P  sin{ \lambda }
 \text {\, on \,} \goth S  \tag5.6.1 \\
\exp(i \lambda \bk H) \bk P \exp(- i \lambda \bk H) =
-  \bk Q    sin{ \lambda}  +  \bk P   cos{ \lambda }
 \text{\,  on \, } \goth S  \tag5.6.2 \\
\exp(i \lambda \bk H) \bk A \exp( -i \lambda \bk H ) =
\exp (-i \lambda) \bk A \text{\,  on  \,} \goth S 
 \tag5.6.3 \\
\exp(i \lambda \bk H) \exp( i \nu \bk Q) \exp(-i \lambda \bk H) =
 \exp\bigl( i \nu ( \bk Q  cos{ \lambda}  +   \bk P  sin{ \lambda
}) \bigr)
 \text {\, on \,}  \frak H  \tag5.6.4 \\
\exp(i \lambda \bk H) \exp(i \nu \bk P) \exp(- i \lambda \bk H)  =
\exp \bigl(i \nu  (- \bk Q sin{ \lambda} + \bk P cos{ \lambda } )
\bigr)
\text{\, on \,}   \frak H  \tag5.6.5
\endgather $$
\endproclaim
\demo{Proof}  Let $ F(\lambda)$ and $ G(\lambda) $ equal the left
sides
of (5.6.1) and (5.6.2). Then
$$ \align
&F^\prime(\lambda) =
 -i \exp(i \lambda \bk H) [\bk Q , \bk H ] \exp(- i \lambda \bk H)
=
G(\lambda) \text {  and } \\
&G^\prime(\lambda) =
 -i \exp(i \lambda \bk H) [\bk P , \bk H ] \exp(- i \lambda \bk H)
=
-F(\lambda).
\endalign $$  Also $ F(0) = \bk Q $ and $ G(0) = \bk P $.
One obtains (5.6.1) and (5.6.2) by solving the system of
differential equations for F and G on $\goth S $. (5.6.3) follows
from (5.6.1) and (5.6.2).
\endgraf  Since $\goth S $ is a dense set of analytic vectors for
$\bk P, \bk Q, \bk H$ one deduces from (5.6.1) and (5.6.2)
that  (5.6.4) and (5.6.5) hold
when acting on a fixed vector in $\goth S$ provided $ \vert \nu
\vert
$ is small. Both sides are groups of unitary operators so (5.6.4)
and (5.6.5) hold generally.
\enddemo
 From now on we assume that the Bose-like oscillator is
irreducible.
\proclaim{\jump 5.7 Structure  Theorem for the Bose-like
Oscillator, I}
 \newline        Suppose that
 $(\frak H, \bk P, \bk Q, \bk H ) $ is an irreducible Bose-like
oscillator. Then there exists a real number
 $ \mu \in (-\half , \infty) $ and $ \phi_0 \in \frak H $
as follows :
$$  \bk H \text{  has pure point spectra. The spectrum  } \sigma (
\bk H)
\text { of } \bk H \text { is given by  }   $$
$$ \sigma(\bk H)  =  \{ \mu + \half, \mu + \tfrac32, \mu +
\tfrac52, \dots \}
 .\tag5.7.1 $$
$$\align &\mu + \half \text{ is the smallest eigenvalue of  } \bk
H .
\text { Select } {\phi}_0 \in \frak H \text{ so } \bk H \phi_0 =
(\mu + \half) \phi_0 , \\
& \Vert \phi_0 \Vert  = 1 . \endalign $$
$$ \text{  Define }  \phi_n :=
   ( {\gamma}_\mu(n))^{- 1/2} {\bk A}^{\star n}{\phi}_0, \,
n \in \Bbb N \text { and } {\phi}_{-1} = 0 .
\tag5.7.2 $$
$$ \{ \phi_n \}_{n \in \Bbb N} \text { is a complete orthonormal
set in }
\goth H . \tag 5.7.3 $$
(5.7.4) - (5.7.9) hold  for all $ n \in \Bbb N $.
$$ \align  {\bk A}^{m} {\bk A}^{\star n} \phi_0 &=
   \frac { \gamma_\mu(n)}{ \gamma_\mu(n-m)} 
         {\bk A}^{\star (n - m)} \phi_0  \,
        \text{ if} \, \,  n \ge m  \ge 0             \tag 5.7.4 \\
         &= 0 \text { if } \, m > n         \\
    {\bk A}^\star {\bk A} \phi_n =& ( n + 2 \mu \theta_{n}) \phi_n
\\
    {\bk A}{\bk A}^\star \phi_n =&  (n + 1 + 2 \mu \theta_{n+1} )
\phi_n \\
     \bk H \phi_n =& \half ( \bk A {\bk A}^\star +
         {\bk A}^\star \bk A ) \phi_n =
 (n + \mu + \half) \phi_n . \tag5.7.5
\endalign $$
  {  \it  The commutation relation for the Bose-like oscillator:}
$$ i (\bk P \bk Q - \bk Q \bk P ) =  \bk I  + 2 \mu \bk J  \text {
on  }
\goth S , \text {  where } \bk J :=
 \exp \bigl( - \pi i ( \bk H - ( \mu + \half) \bk I ) \bigr) .
\tag5.7.6 $$
$$  \bk J \bk P = - \bk P \bk J, \, \text { and} \,  \bk J \bk Q =
 - \bk Q \bk J \, \text{ on } \, \goth S \tag 5.7.7 $$ 
$$ \bk J \phi_n = (-1)^n \phi_n , \,  n \in \Bbb N .\tag5.7.8 $$
$$ \bk J = {\bk J}^\star = {\bk J }^{-1} ,  \text  {  and }
\exp \bigl( -2 \pi i ( \bk H - \mu - \half) \bigr) = \bk I .
  \tag5.7.9 $$

\endproclaim
\demo{Proof} Consider the operator $ {\bk J}_0 := \exp(-i \pi \bk
H) $
and assume that $ \nu \in R $ .  Then from (5.6.4) and
 (5.6.5) we get
$ {\bk J}_0^\star \exp(i \nu \bk Q) {\bk J}_0 = \exp(- i \nu \bk Q)
$
     and 
$ {\bk J}_0^\star \exp(i \nu \bk P) {\bk J}_0
\allowmathbreak  = \exp(- i \nu \bk P) . $
Thus $ \bk J_o \bk Q = -\bk Q \bk J_o \,  \text{ and } \, \bk J_o
\bk P
=- \bk P \bk J_o \, \text { on} \, \goth S . $ 
Also 
$ \exp( i \nu \bk Q) { \bk J}_0^2  =$
  ${ \bk J}_0^2 \exp( i \nu \bk Q) $    and
$ \exp( i \nu \bk P) { \bk J}_0^2  =  { \bk J}_0^2   \exp( i \nu
\bk P). $
 It follows from definition 5.1 iv) that $ {\bk J}_0^2 =  c \bk I
$
 for some complex number c. But since  ${\bk J}_0^2 $ is a unitary
operator
 necessarily ${\bk J}_0^2  = \exp(-2 \pi i \alpha) \bk I $ for some
real number
 $ \alpha.$ Next set $\bk J = \exp( \pi i \alpha) {\bk J}_0 =
 \exp (-i \pi (\bk H - \alpha \bk I)). $  Clearly (5.7.7) is true. 
 Since $ {\bk J}^2 = 
\exp (-2 \pi i ( \bk H - \alpha \bk I )) = \bk I,$  we deduce that
(5.7.9) holds provided $ \mu + \half := \alpha $ . When applied
to (5.7.9) the spectral mapping theorem implies that
$ \sigma (\bk H) \subseteq { \mu + \half + \Bbb Z } ,$ where $\Bbb
Z $
is the set of integers. Now, $ \bk H $ is a non-negative operator,
so  we infer that (5.7.1) is valid if  $\lq = \rq $ is replaced by
$\lq \subseteq \rq $. Clearly $ \mu + \half \ge 0 $.
\endgraf Select $ \phi_0 $ and $ \phi_n $ as in the statement of
the
theorem. From (5.4.2) we see that $ \bk A \bk H \phi_0 =
 (\bk H + \bk I) \bk A \phi_0 $
, so $ (\bk H -(\mu - \half)\bk I) \bk A \phi_0 = 0 $. But $ \mu +
\half $
is the smallest eigenvalue of $\bk H$, so $\bk A \phi_0 = 0 $. 
Also
$ [\bk A, \bk A^\star] \phi_0  = \bk A \bk A^\star \phi_0 = 2 \bk
H \phi_0
= ( 2 \mu + 1) \phi_0. $ Notice that $ \bk H \phi_0 =( \mu + \half)
\phi_0$
implies that  $ \bk J \phi_0 = \phi_0 .$ \endgraf
From (5.5.2) and (5.5.3) we have
$ \bk A {\bk A}^{ \star n} \phi_0 - {\bk A}^{\star n} \bk A \phi_0
=
(n + 2 \mu \theta_n) { \bk A}^{\star (n -1)} \phi_0 $ for all
 $ n \in \Bbb N \setminus \{0 \}, $  so by (2.2.7)
 $ \bk A {\bk A}^{ \star n} \phi_0 =
\frac { \gamma_\mu(n) }{ \gamma_\mu(n -1) } {\bk A}^{\star (n - 1)}
\phi_0. $
This in turn implies the first statements in  (5.7.4), and the rest
of
(5.7.4) follows from this and (2.2.7). 
(5.7.5) is a consequence of (5.7.4). \endgraf
Assume $ m \le n ,\, m,n \in \Bbb N .$ Then, using (5.7.4), the
inner product
$ {\langle \phi_n , \phi_m \rangle} =  
{ \bigl( \gamma_\mu(n) \gamma_\mu(m) \bigr)^{-1/2} }
 {\langle  {\bk A}^{\star m} {\bk A}^{\star n} \phi_0 , \phi_0
\rangle}
= \delta_{m,n}. $ Thus
$\{ \phi_n \}$ is an orthonormal set in $\frak H $. Let $\frak M$
be the
linear span of $ {\bk A}^{ \star n }\phi_0, n \in \Bbb N.$ Then
(5.7.4)
implies that $ \bk A $ and ${ \bk A}^\star $ map $\frak M$ into
itself. Thus
$\bk Q$ and $\bk P$ map $\frak M$ into itself, and since the
oscillator
is irreducible, necessarily the closure of $ \frak M$ equals $
\frak H .$ 
Hence (5.7.3) is proved.               
\endgraf  Finally
$  i (\bk P \bk Q - \bk Q \bk P ) \phi_n  
=  ( \bk A {\bk A}^\star -{\bk A}^\star \bk A ) \phi_n = 
( 1 + 2 \mu  \theta_{n+1} - 2 \mu \theta_n ) \phi _n = 
( \bk I  + 2 \mu \bk J ) \phi_n  $
for all $ n \in \Bbb N ,$
proving the commutator identity in (5.7.6)
\endgraf The non-negativity of $\bk H $ implied that $ \mu + \half
\ge 0 $.
Suppose that $ \mu + \half = 0 $ . Then by (5.7.4) $\bk A \phi_0 =
{\bk A}^\star \phi_0 = 0 $. It follows that $\bk P$ and $\bk Q $
both commute with the projection on $ \phi_0 $, contradicting  the
irreducibility assumption of Definition 5.1iv) . It follows that
necessarily $ \mu \in (-\half , \infty ).$
\enddemo
\proclaim{\jump 5.8 Lemma} Suppose $n \in \Bbb N \setminus \{ 0 \}
.$
$$ \gather
 2^{-\half} {\bk A}^{\star} \phi_0 = {\bk Q} \phi_0 = -i \bk P
\phi_0
\text {\quad and \quad  } i[\bk P, \bk Q]\phi_0 =
 ( 1 + 2\mu) \phi_0 \tag5.8.1 \\
i[\bk P , {\bk Q}^n ] \phi_0 = \frac{\gamma _\mu (n)}{ \gamma _\mu
(n-1)}
{\bk Q}^{n-1} \phi_0 \quad , \quad 
 i[{\bk P}^n , \bk Q ] \phi_0 = \frac{\gamma _\mu (n)}{\gamma _\mu
(n-1)}
{\bk P}^{n-1} \phi_0 \tag5.8.2  \\
  [\bk A , {\bk A}^{\star n} ] \phi_0 = \frac{\gamma _\mu (n)}
{\gamma _\mu (n-1)} {\bk A}^{\star( n-1) } \phi_0 \tag5.8.3
\endgather $$
\endproclaim
\demo {Proof} (5.8.1) is true since $ (\bk Q + i \bk P)\phi_0 =
2^\half \bk A \phi_0 = 0 ,$ and $  \bk Q - i \bk P = 2^\half \bk
A^{\star}.$
From (5.5.4) and (5.7.6) we obtain $ i(\bk P \bk Q^n - \bk Q^n \bk
P)\phi_0
=  \bigl( n \bk Q^{n - 1} + 2 \mu \theta_n \bk Q^{ n - 1} \bk J
\bigr)
\phi_0 .$ By (5.7.8) this yields the first equation in (5.8.2).
The other equations in (5.8.1),
  (5.8.2) and (5.8.3) have similar proofs.
\enddemo

\proclaim{\jump 5.9 Theorem} Suppose $p(\cdot) $ is a complex
polynomial
and $ \goth D_\mu $ is the generalized differentiation operator of
(2.4.1).
Then
$$ 
 i [ \bk P , p(\bk Q) ] \phi_0 = ( \goth D_\mu p)(\bk Q) \phi_0
,\quad  
 i [p(\bk P) , \bk Q  ] \phi_0 = ( \goth D_\mu p)(\bk P) \phi_0  $$
$$ [\bk A , p(\bk A^\star)] \phi_0 = (\goth D_\mu p)(\bk
A^\star)\phi_0
 $$ In case $ p(x) = H_n^\mu( \lambda x)$ one has
$ ({\goth D}_{\mu} p)(x) = 2 \lambda n H_{n - 1}^\mu(\lambda x), 
\, \lambda \in \Cal C . $
\endproclaim
\demo{ Proof} By (2.5.2)  $ p(x) = \sum c_n x^n$ is mapped
by $\goth D_\mu$ to $ \sum \frac { \gamma_\mu(n)}{ \gamma_\mu(n -
1)}
c_n x^{n - 1} $, so the theorem follows easily from lemma 5.8.
The last statement is proved in (2.6.1).
\enddemo
\proclaim{\jump 5.10 Theorem} Suppose $ n \in \Bbb N. $  Then the 
following formulas hold:
{  \it  Rodrigues formula:}
$$  {\bk P}^n \phi_0 = i^n \frac { \gamma_\mu(n)}{2^{n/2} n!}
 H_n^\mu(2^{-\half} \bk Q) \phi_0  
= i^n \frac{\gamma_\mu(n)}{2^n n!} H_n^\mu(2^{-\half}
 {\bk A}^{\star }) \phi_0 \tag5.10.1 $$
{ \it  Dual Rodrigues formula:} 
 $$ 
 {\bk Q}^n \phi_0 =(- i)^n \frac { \gamma_\mu(n)}{2^{ n/2} n!}
 H_n^\mu(2^{-\half} \bk P) \phi_0  =
 (-i)^n \frac{ \gamma_\mu(n)}{ 2^n n!}
 H_n^\mu(i 2^{-\half}  {\bk A}^{\star}) \phi_0 \tag 5.10.2 $$
$$ {\bk A}^{\star n} \phi_0 = \frac { \gamma_\mu(n)}{2^{n/2}n!}
 H_n^\mu(\bk Q) \phi_0 \tag 5.10.3 $$
$$
{\phi_n} =  \frac{ (\gamma_\mu(n))^{\half}}{2^{n/2}n!}
H_n^\mu(\bk Q)   \phi_0  
         = (-i)^n  \frac{ (\gamma_\mu(n))^{\half}}{2^{n/2}n!}
H_n^\mu(\bk P)   \phi_0 \tag 5.10.4
 $$
\endproclaim
\demo{ Proof}We prove (5.10.1) using Theorem 5.9 and the three term
recursion relation (2.6.3) written in the form
$$ 2^{-\half} \frac{\gamma_\mu(n+1)}{(n+1)\gamma_\mu(n)
}H_{n+1}^\mu
(2^{-\half}x) = x H_n^\mu(2^{-\half} x) 
 - 2^{\half} n H_{n-1}^\mu(2^{-\half} x).\text { \quad Then } $$
$ i [ \bk P, H_n^\mu(2^{-\half}\bk Q) ] \phi_0 =
2^{\half} n H_{n-1}^\mu(2^{-\half}\bk Q) \phi_0, $ \,  so \,
$i \bk P H_n^\mu(2^{-\half} \bk Q) \phi_0 +
  \bk Q  H_n^\mu(2^{-\half}\bk Q) \phi_0  \allowmathbreak =
 2^{\half} n H_{n - 1}^\mu(2^{-\half} \bk Q) \phi_0 .$
 The induction proof of (5.10.1) proceeds from this
equation. \endgraf
The remaining equations have similar proofs using 5.9 and (2.6.2).
\enddemo
\proclaim{\jump 5.11 Structure  Theorem for the Bose-like
Oscillator, II}
 Suppose that
 $(\goth H, \bk P, \bk Q, \bk H ) $ is an irreducible Bose-like
oscillator, and  $ \mu $ and $ \phi_n $
are as in Theorem 5.7. Define the unitary operator
$ \Cal F $ on $ \goth H $ to $ \goth H $ by
$$ \align
  \Cal F   &=
\exp \bigl( -\half \pi i ( \bk H - ( \mu + \half ) \bk I ) \bigr)
.
\text{\quad  Then }   \tag5.11.1    \\
{\Cal F}^2   &=   \bk J \text { \quad and } \quad { \Cal F}^\star
= 
\bk J {\Cal F} = {\Cal F} \bk J  ; \tag5.11.2      \\
\bk P &= {\Cal F}^\star \bk Q \Cal F \text {  on } \goth S ;
\tag5.11.3  \\
\Cal F \phi_n &= (-i)^n \phi_n ,\,  n \in \Bbb N . \tag5.11.4
\endalign  $$
\endproclaim
\demo{ Proof } Compare the definitions of $\Cal F$ and $\bk J$
in (5.11.1)  and (5.7.6) to derive (5.11.2). (5.11.3) comes
from (5.6.1) with $\lambda = \pi/2 .$  The eigenvectors of
 $\bk H$ are necessarily the eigenvectors of $\Cal F$ and
thus (5.11.4) is true. 
\enddemo
We show, finally, that given any abstract irreducible Bose-like
oscillator, there exists a number $\mu \in (-\half, \infty)$
such that the abstract Bose-like oscillator is unitarily
equivalent to the concrete irreducible
Bose-like oscillator on  $L_\mu^2(R)$ specified in 3.4 to 3.7.
Thus one has a generalization of the von Neumann uniqueness
theorem \cite{26,p275}. Formal aspects of the  physical theory
are detailed in Ohnuki and Kamefuchi \cite {24}, Chapter 23,
entitled \it { The wave-mechanical representation for a 
Bose-like oscillator }. 
\proclaim{\jump 5.12 Representation Theorem for the Bose-like
Oscillator } Suppose
 $(\goth H, \bk P, \bk Q, \bk H ) $ is an irreducible Bose-like
oscillator, and  maintain the notation of Theorems 5.7 and 5.11.
Define the unitary mapping $ \bk U $ of $ \goth H $ onto the
Hilbert space $  L_\mu^2(R)$  by $ \bk U :
\phi_n  \longmapsto  \phi_n^\mu , n \in \Bbb N $ .
 Then $\bk U$ maps $ \bk P , \bk Q, \bk H,
\bk J,    \text{  and  } \Cal F $ onto  
 $ \bk P _\mu , \bk Q _\mu ,  \bk H_\mu , \bk J_\mu,
 \text {  and  } \Cal F _\mu $ respectively .
\endproclaim
\demo{ Proof}
(5.7.3) and (3.5.2) assure us that $\bk U $ maps a complete
orthonormal set onto a complete orthonormal set and thus $\bk U $
is a unitary mapping. In addition,  $\bk U $
maps $ \bk A $ of 5.2   onto $ \bk A_\mu $ of (3.4.3) because
of the action of  these operators on the orthonormal sets, see
(5.7.2) and (3.7). The other assertions then follow.
\enddemo

\enddocument
\bye